\documentclass{amsart}

\usepackage{amssymb}

\renewcommand{\leadsto}{\mapsto} 
\newcommand{\MU}[1]{m_{[#1]}}
\newcommand{\sr}{{\operatorname{\mathbf {s}}}}
\newcommand{\bC}{{\operatorname{\mathbf{ C}}}}

\newcommand{\ZFCa}{{\operatorname{\mathsf {ZFC}}}}
\newcommand{\CH}{\operatorname{\mathsf {CH}}}
\newcommand{\BC}{\operatorname{\mathsf {BC}}}
\newcommand{\DBC}{\operatorname{\mathsf {DBC}}}

\newcommand{\perf}{\operatorname{\mathsf {Perf}}}
\newcommand{\bv}[1]{[\![#1]\!]}
\newcommand{\wh}[1]{\overline{#1}}

\newcommand{\reals}{{\mathbb R}}
\newcommand{\rationals}{{\mathbb Q}}
\newcommand{\rest}{{\mathord{\restriction}}}

\newcommand{\cov}{\operatorname{\mathsf  {cov}}}

\newcommand{\COV}{\operatorname{\mathsf   {COV}}}

\newcommand{\FF}{\operatorname{\mathsf   {F}}}
\newcommand{\bF}{\operatorname{\mathbf   {F}}}
\newcommand{\bG}{\operatorname{\mathbf   {G}}}

\newcommand{\dom}{{\operatorname{\mathsf {dom}}}}

\newcommand{\suc}{{\operatorname{\mathsf    {succ}}}}

\newcommand{\QED}{\hspace{0.1in} \square \vspace{0.1in}}

\newcommand{\N}{{\mathcal N}}
\newcommand{\M}{{\mathcal M}}

\renewcommand{\>}{\rangle}

\newcommand{\Proof}{{\sc Proof}. \hspace{0.2in}}

\newcommand{\lft}[2]{\mathopen\ifcase#1{}\oo\or
                        \big#2\or\Big#2\else\oo\fi} 
\newcommand{\rgt}[2]{\mathclose\ifcase#1{}\oo\or
                        \big#2\or\Big#2\else\oo\fi}

\renewcommand{\t}{{\mathbf t}}
\newcommand{\p}{{\mathbf p}}
\newcommand{\Ff}{{\mathbf F}}

\newcommand{\q}{{\mathbf q}}

\newcommand{\SM}{{\mathcal{SM}}}
\newcommand{\SN}{{\mathcal  {SN}}}

\theoremstyle{plain}
\newtheorem{theorem}{Theorem}[section]
\theoremstyle{plain}
\newtheorem{lemma}[theorem]{Lemma}

\newtheorem{definition}[theorem]{Definition}

\begin{document}

\title{Strongly meager sets do not form an ideal}
\author{Tomek Bartoszynski}
\address{Department of Mathematics and Computer Science\\
Boise State University\\
Boise, Idaho 83725 U.S.A.}
\thanks{First author  partially supported by 
NSF grant DMS 95-05375 and Alexander von Humboldt Foundation}
\email{tomek@math.idbsu.edu, http://math.idbsu.edu/\char 126 tomek}
\author{Saharon Shelah}
\thanks{Second author partially supported by Basic Research Fund,
Israel Academy of Sciences, publication 607}
\address{Department of Mathematics\\
Hebrew University\\
Jerusalem, Israel}
\email{shelah@sunrise.huji.ac.il, http://math.rutgers.edu/\char 126
  shelah/}

\begin{abstract}
A set $X \subseteq \reals$ is strongly meager if for every measure
zero set $H$,  $X+H \neq \reals$. Let $\SM$ denote the collection of
  strongly meager sets. We show that assuming $\CH$, 
$\SM$ is not an ideal.
\end{abstract}
\maketitle

\section{Introduction}
In 
1919 Borel wrote the paper \cite{Borel} in which he attempted to
classify all measure zero subsets of the real line. 
In this paper he
introduced  a  class of measure zero sets, which are now called strong
measure zero sets.
In 70's Galvin, Mycielski and Solovay found a characterization of
strong measure zero sets that was formulated using only the concept of 
a first category set and of a translation. That allowed, after
replacing first category with measure zero, to define a dual notion of 
a strongly meager set.
It was expected that the global properties of both families of sets
will be similar. Several results listed below support this
expectation. Nevertheless additive properties of both families of sets 
are different. It is well known that the family of  strong measure
zero sets forms an ideal, i.e. is closed under finite unions.
The result  of this paper is that, assuming continuum hypothesis, the
collection of strongly meager sets is not closed under finite unions.

In this paper we  work exclusively in the space $2^\omega $
equipped with 
the standard product measure denoted as $\mu$. Let $\N$ and $\M$ denote the
ideal of all $\mu$--measure zero sets, 
and meager subsets of $2^\omega $, respectively.
For $x,y \in 2^\omega$, $x+y \in 2^\omega $ is
defined as $(x+y)(n) = x(n)+y(n) \pmod 2$. In particular, $(2^\omega , 
\operatorname{+})$ is a group and $\mu$ is an invariant measure.

\begin{definition}
A set $X$ of real numbers or more generally, a metric space, is strong
measure zero  if, for each sequence $\{\varepsilon_n: n \in \omega\}$ of
positive real numbers there is a sequence $\{X_n: n \in \omega\}$ of
subsets of $X$ 
whose union is $X$, and for each $n$ the
diameter of $X_n$ is less than $\varepsilon_n$. 
\end{definition}

The family of strong measure zero
subsets of $2^\omega $ is denoted by $\SN$.

The following characterization of strong measure zero 
 is the starting point for our considerations.
\begin{theorem}[\cite{GMSSmz}]\label{solo}
The following  are equivalent:
\begin{enumerate}
\item $X \in \mathcal  {SN}$,
\item for every set $F \in \M$, $X+F \neq 2^\omega$.~$\QED$
\end{enumerate}
\end{theorem}

This theorem indicates that the notion of strong measure
zero should have its category analog.
Indeed, we define after Prikry:
\begin{definition}\label{defstrmea}
  Suppose that $X \subseteq 2^\omega $.

We say that $X$  is strongly meager if for every $H \in
  \N$, $X+H \neq 2^\omega $. Let $\SM$ denote the collection of
  strongly meager sets.
\end{definition}

Observe that if $z \not\in X+F=\{x+f: x \in X, f\in F\}$ then $X \cap
(F+z) = \emptyset$. In particular, a strong measure zero set can be
covered by a translation of any dense $G_\delta $ set, and every
strongly meager set can be covered by a translation of any measure one 
set.

If $X \subseteq 2^\omega $ is a group then the concepts of strong
measure zero and strongly meager connect to the classical construction of a
nonmeasurable set by Vitali (a selector of $\reals/\rationals$).

\begin{theorem}[Reclaw]
  Suppose that $X \subseteq 2^\omega $ is a dense subgroup of
  $(2^\omega,+)$.
Then
\begin{enumerate}
\item $X \in \SM$ if and only if every selector from $2^\omega / X$ is
  nonmeasurable.
\item $X \in \SN$ if and only if every selector from $2^\omega / X$
  does not have the Baire property.
\end{enumerate}
\end{theorem}
\Proof
The proof below requires the  group $X$ to be infinite and the set $2^\omega
/ X$ to be infinite. A dense group will have these properties.

We will show only (1), the proof of (2) is analogous.
Note that if $X$ is a selector from $2^\omega / X$ and $X$ is as above 
then $X$ is nonmeasurable if and only if $X$ does not have measure
zero. 

$ \rightarrow $ Suppose that $X \in \SM$ and $H \in \N$. 
Let $x \not \in X+H$. It follows that $[x]_X \cap H =\emptyset$, hence
no selector is contained in $H$.

\bigskip

$ \leftarrow$ Suppose that $X \not \in \SM$ and let $H \in \N$ be such
that $X+H=2^\omega $. For each $x \in 2^\omega $, $[x]_X \cap H \neq
\emptyset$. It follows that we can choose a selector contained in
$H$.~$\QED$

Note that $X \not \in \SN$ if there exists a meager set $F$ such that
the family $\{F+x: x \in X\}$ covers $2^\omega $. 
Instead of the assignment $x \leadsto F+x$ we can consider a more
general mapping $x \leadsto (H)_x$, where $H \subseteq 2^\omega \times 
2^\omega $ is a Borel set such that $(H)_x = \{y:\<x,y\>\in H\}\in \M$ 
for all $x \in 2^\omega $.
\begin{definition}
$X \in \COV(\M)$ if
for every Borel set $H \subseteq 2^\omega \times 
2^\omega$ such that $(H)_x \in \M$ for all $x \in 2^\omega $,
$$\bigcup_{x \in X} (H)_x \neq
2^\omega.$$
 Similarly,
$X \in \COV(\N)$ if
for every Borel set $H \subseteq 2^\omega \times 
2^\omega$ such that $(H)_x \in \N$ for all $x \in 2^\omega $,
$$\bigcup_{x \in X} (H)_x \neq
2^\omega.$$

\end{definition}

Note that 
\begin{lemma}
  $\COV(\N) \subseteq \SM$ and   $\COV(\M) \subseteq \SN$.
\end{lemma}
\Proof
Given $F \in \M $ let $H=\{(x,y): y \in F+x\}$.
It is clear that,
$\bigcup_{x \in X} (H)_x = F+X$.~$\QED$

Families $\SN$ and $\SM$ as well as  $\COV(\M)$ and $\COV(\N)$ are
dual to each other and we are interested to what extent the properties 
of one family are shared by the dual one.

Below we present several results of that kind.
The proofs of these results as well as quite a lot of additional
material can be found in \cite{BJbook}.

\begin{definition}
  Let Borel Conjecture ($\BC$) be the assertion that there are no uncountable
  strong measure zero sets, and
Dual Borel Conjecture ($\DBC$) be the assertion that there are no
uncountable
  strongly meager sets.
\end{definition}

Sierpinski showed that Borel Conjecture contradicts $\CH$. His proof essentially
yields the following:
\begin{theorem}
  Assume $ {\bf MA} $. 
Both $\COV(\M)$ and $\COV(\N)$ contain sets of size $
2^{\boldsymbol\aleph_0} $. In particular, both Borel Conjectures are false. 
\end{theorem}

There are many weaker assumptions than $ {\bf MA} $ that contradict
$\BC$ or $\DBC$. Nevertheless we have the following:

\begin{theorem}[\cite{Lav76Con}]
  Borel Conjecture is consistent with $\ZFCa$.
\end{theorem}

\begin{theorem}[\cite{CarStr}]
  Dual Borel Conjecture is consistent with $\ZFCa$.
\end{theorem}

\begin{definition}
  An uncountable set $X \subseteq 2^\omega $ is a Luzin set if $X \cap F$ is
countable for $F \in \M$, and 
 is a Sierpinski set  if $X \cap G$ is
countable for $G \in \N$.
\end{definition}

Sierpinski showed that every Luzin set is in $\SN$.
In addition we have the following:
\begin{theorem}[ \cite{RecOpen}]
Every Luzin set is in $\COV(\M)$.
\end{theorem}

\begin{theorem}[\cite{Paw92Sie}]
  Every Sierpinski set is in $\COV(\N)$ (and so in $\SM$).
\end{theorem}

Results presented above indicate that we have certain degree of
symmetry between the notions of strongly meager and strong measure
zero. 
The main objective of this paper is to show that as far as additive
properties of both families are concerned it is not the case. 

Sierpinski showed that $\SN$ is a $ \sigma $-ideal. 
 In fact, we have the following:
 \begin{theorem}[\cite{CarStr}]
   Assume ${\bf MA} $. Then the additivity of $\SN$ is
  $2^{\boldsymbol\aleph_0}$.
 \end{theorem}
Similarly,
\begin{theorem}[\cite{BarJud93Cov}]
  \begin{enumerate}
  \item $\COV(\M)$ is a $ \sigma $-ideal, 
  \item Assume ${\bf MA} $. Then the additivity of $\COV(\M)$ is
  $2^{\boldsymbol\aleph_0}$.
  \end{enumerate}
\end{theorem}

Surprisingly the dual results are not true.
\begin{theorem}
  It is consistent that $\COV(\N)$ is not a $ \sigma $-ideal.
\end{theorem}
\Proof
It is an immediate consequence of the following theorem 
of Shelah:
\begin{theorem}[\cite{Sh592}]
  It is consistent that
  $\cov(\N)=\boldsymbol\aleph_\omega$.
\end{theorem}
Recall that
$$\cov(\N )=\min\left\{|{\mathcal A}|: {\mathcal A} \subseteq 
{\mathcal N}  \ \&\ \bigcup {\mathcal A} =2^\omega \right\}.$$

Suppose that $\cov(\N)=\boldsymbol\aleph_\omega$ and let a family ${\mathcal A} 
\subseteq \N$ witness that.
Let $H \subseteq 2^\omega \times 2^\omega $ be an Borel set with null
vertical sections and  such that 
$$\forall G \in \N \ \exists x \in 2^\omega \ G \subseteq (H)_x.$$
Such a set can be easily constructed from a universal set.

For each $G \in {\mathcal A} $ choose $x_G \in 2^\omega $ such that 
$G \subseteq (H)_{x_G}$. It follows that  $X=\{x_G : G \in {\mathcal
  A}\} \not\in \COV(\N)$. On the other hand, every set of size
$<\cov(\N)$ belongs to $\COV(\N)$ and $X$ is a countable union of such 
sets.~$\QED$

The purpose of this paper is to show that
\begin{theorem}\label{biggie}
  Assume $\CH$. Then $\SM$ is not an ideal.
\end{theorem}

\section{Framework}\label{out}
The proof of Theorem \ref{biggie}  occupies the rest of the paper. The
construction is motivated by the tools and
methods developed in \cite{RoSh470}. We should note here
that by using the forcing notion defined in this paper we can also 
show that the statement ``$\SM$ is not an
ideal'' is not equivalent to $\CH$.  
However, since the main result is of
interest outside of set theory we present a version of the
proof that does not contain any metamathematical references. 

The structure of the proof is as follows:
\begin{itemize}
\item In section \ref{out} we show 
that in order to show that $\SM$ is not an ideal it suffices to find
certain partial ordering ${\mathcal P} $ (Theorem \ref{newtrick}).
\item  The definition of $ {\mathcal P} $ involves construction of a
  measure zero set $H$ with some special properties. All results
  needed to define $H$ are proved in section \ref{three}, and $H$
  together with other parameters is
  defined in section \ref{four}.
\item ${\mathcal P} $ is defined in section \ref{seven}. The proof
  that $ {\mathcal P} $ has the required properties is a consequence
  of Theorem \ref{crucialgeneralmore}, which is the main result of
  section \ref{five}, and Theorems \ref{main} and \ref{main1}, which
  are proved in section \ref{six}.
\end{itemize}

We will show that in order to prove \ref{biggie} it is enough to
construct a partial ordering satisfying several general
conditions. Here is the first of them.
\begin{definition}
Suppose that $ ({\mathcal P}, \geq)$ is a partial ordering. 
We say that ${\mathcal P} $ has the fusion property if there exists a
sequence of binary relations $\{ \geq_n : n \in \omega\}$ (not
necessarily transitive) such that
\begin{enumerate}
\item If $p \geq_{n} q$ then $p \geq q$, 
\item if $p \geq_{n+1} q$ and $r \geq_{n+1} p$ then $r \geq_n q$,
\item if $\{p_n: n \in \omega\}$ is a sequence such that $p_{n+1}
  \geq_{n+1} p_n$ for each $n$ then there exists $p_\omega $ such that
  $p_\omega \geq_n p_n$ for each $n$.
\end{enumerate}
\end{definition}

From now on we will work in $2^\omega $ with the set of rationals defined as 
$$\rationals=\{x
\in 2^\omega : \forall^\infty n \ x(n)=0\}.$$

Let $\perf$ be the collection of perfect subsets of $2^\omega $.
For $p, q \in \perf $ let $p \geq q$ if $p \subseteq q$.

We will be interested in subsets of $\perf \times \perf$. Elements of
$\perf \times \perf$ will be denoted by boldface letters and if
${\mathbf p} \in \perf \times \perf$ then ${\mathbf p}=(p_1, p_2)$.
Moreover, for ${\mathbf p}, {\mathbf q} \in \perf \times \perf$,
${\mathbf p}\geq {\mathbf q}$ if $p_1 \subseteq q_1$ and $p_2
\subseteq q_2$.

\begin{theorem}\label{newtrick}
Assume $\CH$, fix a measure zero set $H \subseteq 2^\omega $, and
suppose that there exists  
 a family $ {\mathcal P} \subseteq \perf \times \perf$ 
such that:
\begin{enumerate}
\item [(A0)] $ {\mathcal P} $ has the fusion property,
\item[(A1)] For every $ {\mathbf p} \in
  {\mathcal P}$, $n \in \omega $ and $z \in 2^\omega $ there exists
  ${\mathbf q} \geq_n {\mathbf p}$ such that 
$q_1 \subseteq H+z$ or $q_2 \subseteq H+z$,
\item[(A2)] for every ${\mathbf p} \in {\mathcal P} $, $ n \in \omega
  $, $X \in
  [2^\omega]^{\leq \boldsymbol\aleph_0} $, $i=1,2$ and $\t \in \perf$ such that
  $\mu(\t)>0$,
\begin{displaymath}\mu\lft2(\lft1\{z\in 2^\omega : \exists \q \geq_n \p\  
 X \cup (q_i + 
\rationals) 
\subseteq \t+\rationals+z\rgt1\}\rgt2)=1.\end{displaymath}
\end{enumerate}
Then $\SM$
 is not an ideal.
  \end{theorem}
   \Proof
 We intend to build by induction sets $X_1,X_2 \in \SM$ in such a way
that $H$ witnesses that $X_1 \cup X_2$ is not strongly meager, that
is, $(X_1\cup X_2)+H=2^\omega $. 
By induction we will define an $\omega_1$-tree of members of
${\mathcal P} $ and then take the selector from the elements of this
tree. This is a 
refinement of the method
invented by Todorcevic  (see \cite{GalMil84Gam}), who used an
Aronszajn tree of perfect sets to construct a set of reals with some
special properties. More examples can be found in \cite{BarInv}.

For each $ \alpha< \omega_1 $, $\mathfrak
 T_\alpha $ will denote the $\alpha $'th level of an Aronszajn tree of
 elements of $ 
 {\mathcal P} $.  
More precisely, we will define $\suc({\mathbf p},\alpha)
 \subseteq 
 {\mathcal P}$ -- the
 collection of all successors of ${\mathbf p}$ on level $ \alpha $.
We will require that:
\begin{enumerate}
\item $\mathfrak T_0= \{2^\omega \times 2^\omega\} $,
\item $\suc({\mathbf p}, \alpha ) $ is countable  (so levels of the tree are
  countable),
\item if ${\mathbf q} \in \suc({\mathbf p},\alpha)$ then ${\mathbf q}
  \geq {\mathbf p}$,
\item if $\suc({\mathbf p},\alpha)$ is defined then for each $n \in \omega $
  there is ${\mathbf q} \in \suc({\mathbf p},\alpha)$ such that
  ${\mathbf q} \geq_n {\mathbf p}$. 
\end{enumerate}
Note that the tree constructed in this way will be an 
Aronszajn tree since an uncountable branch would produce an
uncountable descending sequence of closed sets. For an arbitrary $
{\mathcal P} $ with fusion property the
conditions above will guarantee that we build an $ \omega_1$-tree with
countable levels. This suffices for the constructions we are
interested in.

Let ${\mathfrak T}=\bigcup_{\alpha<\omega_1} {\mathfrak T}_\alpha $
 where ${\mathfrak T}_\alpha = \suc(2^\omega \times 2^\omega , \alpha)$.
 For each $ {\mathbf p} \in {\mathfrak T}_\alpha $ choose $x_{\mathbf
 p}^1 \in 
 p_1 $ and $x_{\mathbf p}^2 \in p_2 $. We will show
 that we can arrange this construction in such a way that
 $X_1=\{x^1_{\mathbf p}: {\mathbf p}\in \mathfrak T\}$ and
 $X_2=\{x^2_{\mathbf p}: {\mathbf p}\in \mathfrak T\}$
 are the sets we are looking for.

Let $\{(\t_\alpha,i_\alpha)  : \alpha < \omega_1\}$ be an enumeration
of pairs $(\t,i)\in \perf \times \{1,2\}$
such that $\mu(\t)>0$. Let $\{z_\alpha: \alpha <
 \omega_1\}$ be an enumeration of $2^\omega $.

\bigskip

{\sc Successor step}.

Suppose that $\mathfrak T_\alpha $ is already constructed.
 Denote $X^\alpha =
 \left\{x^1_{\mathbf p},x^2_{\mathbf p}: {\mathbf p} \in \bigcup_{\beta
   \leq \alpha } \mathfrak T_\beta\right\}$.

For each $\p \in {\mathfrak T}_\alpha $ and $ n \in \omega $, let 
$$Z^n_{\p}=\left\{z\in 2^\omega : \exists \q \geq_n \p \ X^\alpha 
\cup  (q_{i_\alpha} + \rationals) 
\subseteq \t_\alpha +\rationals+z\right\}.$$
Note that by A2, each set $Z^n_{\p}$ has measure one.
Fix
$$y_\alpha \in \bigcap_{\p\in {\mathfrak T}_\alpha} \bigcap_{n \in
  \omega} Z^n_{\p}.$$
For each $\p \in {\mathfrak T}_\alpha $ choose $\{{\mathbf p}^n: n \in \omega\}$
such that 
\begin{enumerate}
\item ${\mathbf p}^n \geq_{n+1} {\mathbf p}$ for each $n$,
\item $X^\alpha  \cup
  (p^n_{i_\alpha} + \rationals)  
\subseteq \t_\alpha +\rationals+y_\alpha$.
\end{enumerate}

Next apply A1 to get sets $\{{\mathbf q}^n: n \in \omega\}$ such that for all $n$,
\begin{enumerate}
\item ${\mathbf q}^n \geq_{n+1} {\mathbf p}^n$,
\item $q^n_1 \subseteq H+z_\alpha $ or $q^n_2 \subseteq H+z_\alpha $.
\end{enumerate}
Define $\suc({\mathbf p},\alpha+1)=\{{\mathbf q}^n: n \in \omega\}$.
Note that for each $n \in \omega $
  there is ${\mathbf q} \in \suc({\mathbf p},\alpha)$ such that
  ${\mathbf q} \geq_n {\mathbf p}$. 
For completeness, if $\mathbf p \in \bigcup_{\beta<\alpha } \mathfrak
T_\beta $ then
put 
$$\suc({\mathbf p}, \alpha+1)=\bigcup \{\suc({\mathbf q},
\alpha+1): {\mathbf q} \in \suc({\mathbf p}, \alpha)\}.$$

\bigskip

{\sc Limit step}.

Suppose that $ \alpha $ is a limit ordinal and $\mathfrak T_\beta $
 are already constructed for $ \beta < \alpha $.
Suppose that ${\mathbf p}_0 \in \mathfrak T_{\alpha_0}$,
$\alpha_0<\alpha$. 
Find an increasing sequence   $\{\alpha_n: n \in \omega\}$
with $\sup_n \alpha_n=\alpha $, and for $k \in \omega $,
let $\{{\mathbf p}_n^k : n
\in \omega\}$ be such that 
\begin{enumerate}
\item ${\mathbf p}_n^k \in \mathfrak T_{\alpha_n}$,
\item ${\mathbf p}_{n+1}^k \geq_{n+k+1} {\mathbf p}_n^k$ for each $k,n\in
  \omega $.
\end{enumerate}

Let ${\mathbf p}_\omega^k$ be such that ${\mathbf p}_\omega^k
\geq_{n+k} {\mathbf p}_n^k$.
Define $\suc({\mathbf p}_0,\alpha)=\{{\mathbf p}^k_\omega : k \in
\omega\}$. 
This concludes the construction of $\mathfrak T$ and $X_1,X_2$. 

\bigskip

\begin{lemma}
  $X_1, X_2 \in \SM$.
\end{lemma}
\begin{proof}
  We will show that $X_1 \in \SM$. The proof that $X_2 \in
\SM$ is the same. 

Let $G \subseteq 2^\omega $ be a
  measure zero set. Find $ \alpha < \omega_1$ such that $G \cap (\t_\alpha +
  \rationals)=\emptyset$ and $i_\alpha=1$. 
It follows that,
\begin{displaymath}X_1 \subseteq X^\alpha \cup \bigcup_{{\mathbf p}
    \in {\mathfrak T}_{\alpha+1} } 
p_1 \subseteq \t_\alpha +\rationals+y_\alpha \subseteq (2^\omega
\setminus G)+ y_\alpha .\end{displaymath}
Thus $X_1 + y_\alpha \subseteq 2^\omega
\setminus G$ and therefore $ y_\alpha  \not \in X_1+G$, which finishes the proof.~$\QED$ 

\begin{lemma}
  $X_1 \cup X_2 \not\in \SM$.
\end{lemma}
\Proof
  Let $H$ be the set used in A1.
We will show that $(X_1 \cup X_2) + H=2^\omega $.
Suppose that $z \in 2^\omega $ and let $\alpha<\omega_1$ be such that
$z=z_\alpha $.
By our construction, for any $\mathbf p \in \mathfrak T_{\alpha+1} $, $x^1_{\mathbf p}
\in z +H$ or $x^2_{\mathbf p} \in z
+H$.
Thus $z \in (X_1\cup X_2)+H$, which ends the proof.
\end{proof}

This shows that the sets $X_1, X_2$ and $H$ have the required
properties. The proof of \ref{newtrick} is finished.~$\QED$

Therefore the problem of showing that $\SM$ is not an ideal reduces to the
construction of an appropriate set $ {\mathcal P} $. 
We will do that in the following sections.

\section{Measure zero set}\label{three}
In this section we will develop tools to define a measure zero set $H$ that will
be used in the construction of $ {\mathcal P} $ and will 
witness that the union of two strongly meager sets $X_1, X_2$ defined
 in the proof of \ref{newtrick} is not strongly meager. 
The set $H$ will be defined at the end
of the next section.

We will need several definitions.

\begin{definition}
Suppose that $I \subseteq \omega $ is a finite set. 
Let $\FF^I$ be the collection of all functions $f: \dom(f)
\longrightarrow 2$, with $\dom(f) \subseteq 
2^I$. For $f \in \FF^I$, let $m_f^0=|\{s: f(s)=0\}|$ and $m_f^1=|\{s:
f(s)=1\}|$.

For a set $B \subseteq 2^I$ let $(B)^1 = 2^I \setminus B$ and $(B)^0=B$.

We will work in
the space $(2^I,+)$ with addition mod $2$. For a function $f \in
\FF^I$ let 
\begin{displaymath}(B)^f = \bigcap_{s \in \dom(f)} (B+s)^{f(s)}.\end{displaymath}
In addition let $(B)^\emptyset = 2^I$.

For $f\in \FF^I$ and $ k \in \omega $, let 
\begin{displaymath}\FF^I_{f,k}=\left\{g\in \FF^I: f \subseteq g \ \&\ |\dom(g) \setminus
\dom(f)| \leq k\right\}.\end{displaymath} 
  
\end{definition}
The set $H$ will be defined using an infinite sequence of finite
sets. The following theorem describes how to construct one term of this 
sequence.
\begin{theorem}\label{choosec}
   Suppose that $m \in \omega $ and $0<\delta < \varepsilon < 1$
  are given. There exists $n \in \omega $ such that for every finite
  set $I \in [\omega]^{>n}$ there exists a set $C \subseteq 2^I$ such
  that $1-\varepsilon + \delta \geq |C|\cdot 2^{-|I|} \geq 1-
  \varepsilon -\delta$ and for every $f \in 
  \FF^I_{\emptyset, m}$,
$$\left| \frac{|(C)^f|}{|(C)^\emptyset|} - (1-\varepsilon)^{m^0_f}
  \varepsilon^{m^1_f}\right| < \delta .$$
\end{theorem}
Note that 
the theorem says that we can choose $C$ is such a way that for any
    sequences 
    $s_1, \dots, s_{m} \in 2^I$ the sets $s_1+C,  \dots,
    s_{m}+C$ are probabilistically independent with error $\delta $.
 Thus, we want $ \delta $ to be much smaller than $ \varepsilon^m$.
In order to prove this theorem it is enough to verify the following:

\begin{theorem}\label{choosec1}
  Suppose that $m \in \omega $ and $0<  \delta < \varepsilon < 1$
  are given. There exists $n \in \omega $ such that for every finite
  set $I \in [\omega]^{>n}$ there exists a set $C \subseteq 2^I$ such
  that $1-\varepsilon + \delta \geq |C|\cdot 2^{-|I|} \geq 1-
  \varepsilon - \delta $ and for every set $X 
  \subseteq 2^I$, $|X| \leq m$
\begin{displaymath}\left| \frac{|\bigcap_{s\in X} (C+s)|}{2^{|I|}} -
    (1-\varepsilon)^{|X|}\right| <
    \delta.\end{displaymath}
\end{theorem}
\Proof Note first that \ref{choosec1} suffices to prove
\ref{choosec}. Indeed, if for every $X \in [2^I]^{\leq m}$,
\begin{displaymath}\left| \frac{|\bigcap_{s\in X} (C+s)|}{2^{|I|}} -
    (1-\varepsilon)^{|X|}\right| <
    \delta\end{displaymath} 
then we show by induction on $m^1_f$ that for every $f \in \FF^I_{\emptyset, m}$,
$$\left| \frac{|(C)^f|}{|(C)^\emptyset|} - (1-\varepsilon)^{m^0_f}
  \varepsilon^{m^1_f}\right| < 2^m\delta. $$

\bigskip 

Fix $ m, \delta$ and $ \varepsilon$,   and choose  the set $C
\subseteq 2^I$ randomly (for the moment $I$ is arbitrary). For each $s
\in 2^I$ decisions whether $s \in  
 C$ are made independently with  the probability of $s \in C$ equal to $1-
\varepsilon $. 
Thus the set $C$ is a result of a sequence of Bernoulli trials.
Note that by the Chebyshev's inequality, 
 the probability that $1-\varepsilon + \delta \geq |C|\cdot 2^{-|I|} \geq 1-
  \varepsilon -\delta $ approaches $1$ as $|I|$ goes to infinity. 

Let $S_n$ be the number of successes in $n$ independent Bernoulli
  trials with probability of success $p$.
We will need the following well--known fact that we will prove here
for completeness.
\begin{theorem}\label{cor1}
  For every $ \delta  >0$, 
\begin{displaymath}P\left(\left|\frac{S_n}{n}-p\right| \geq \delta \right) \leq
2 e^{-n\delta^2/4}.\end{displaymath}
\end{theorem}
\Proof
We will show that
\begin{displaymath}P\left(\frac{S_n}{n} \geq p+\delta \right) \leq
e^{-n\delta^2/4}.\end{displaymath}
The proof that
\begin{displaymath}P\left(\frac{S_n}{n} \leq p-\delta \right) \leq
e^{-n \delta^2/4}\end{displaymath}
is the same.
Let $q = 1-p$. Then for each $x \geq 0$ we have 
\begin{multline*}
  P\left(\frac{S_n}{n} \geq p+\delta \right) \leq \sum_{k \geq
    n(p+\delta)}^n {n \choose k} p^k q^{n-k} \leq \\
\sum_{k \geq
    n(p+\delta)}^n e^{-x(n(p+\delta) -k)}\cdot {n \choose k}
  p^k q^{n-k} \leq \\
e^{-xn\delta} \cdot \sum_{k \geq
    n(p+\delta)} {n \choose k} (p e^{xq})^k (q e^{-xp})^{n-k} \leq\\
e^{-xn\delta} \cdot \sum_{k=0}^n {n \choose k} (p e^{xq})^k (q
e^{-xp})^{n-k} =
e^{-xn\delta} \left(p e^{xq} + q e^{-xp}\right)^n \leq \\
e^{-xn\delta} \left(p
e^{x^2q^2} + q e^{x^2p^2}\right)^n \leq 
e^{-xn\delta} \left(p e^{x^2} + q e^{x^2}\right)^n = 
e^{-xn\delta} e^{nx^2} = e^{n(x^2 - \delta x)}.
\end{multline*}
The inequality $p e^{xq} + q e^{-xp} \leq p
e^{x^2q^2} + q e^{x^2p^2}$ follows from the fact that $e^x \leq
e^{x^2}+x$, for every $x$.
The expression $e^{n(x^2 - \delta x)}$ attains its minimal value at
$x = \delta/2$, which yields the desired inequality.~$\QED$

  Consider an arbitrary set $X \subseteq 2^I$. To simplify
  the notation denote $V = 2^I \setminus C$ and note that
$\bigcap_{s \in X} (C+s) = 2^I \setminus (V+X)$.
For a point $t \in 2^I$, $t \not\in X+V$ is equivalent to $(t+X)\cap
V=\emptyset$.
Thus  the probability that $t \not\in X+V$ is equal
to $ (1-\varepsilon)^{|X|}$.

Let $G(X)$ be a
subgroup of $(2^I, +)$ generated by $X$. Since every element of $2^I$
has order $2$, it follows that $|G(X)|\leq 2^{|X|}$.

\begin{lemma}
  There are sets $\left\{U_j : j \leq {|G(X)|}\right\}$ such that: 
  \begin{enumerate}
  \item $ \forall j \ \forall s, t \in U_j \ \lft1(s \neq t
    \rightarrow s+t \not\in G(X)\rgt1)$, 
  \item $ \forall j \leq {|G(X)|} \ |U_j|=2^{|I|}/|G(X)|$,
  \item $ \forall i \neq j \ U_i \cap U_j = \emptyset$,
  \item $\bigcup_{j \leq |G(X)|} U_j = 2^I$.
  \end{enumerate}
\end{lemma}
\Proof
Choose  $U_j$'s to be disjoint selectors from the cosets
$2^I/G(X)$.~$\QED$ 

Note that if $t_1, t_2 \in U_j$ then the events $t_1 \in X+V$ and $t_2
\in X+V$ are independent since sets $t_1+X$ and $t_2+X$ are disjoint.
Consider the sets $X_j=
U_j \cap \bigcap_{s \in X} (C+s)$ for $j
\leq |G(X)|$.
The expected value of the size of this set is $ (1-\varepsilon)^{|X|} \cdot
2^{|I|}/|G(X)|$. 
By \ref{cor1} for each $j \leq |G(X)|$, 
\begin{displaymath}P\left(\left|\frac{|X_j|}{2^{|I|}/|G(X)|} - (1-\varepsilon)^{|X|}\right|
\geq \delta \right) 
  \leq 2 
  e^{-2^{|I|-2}\delta^2/|G(X)|}.\end{displaymath}
It follows that for every $X \subseteq 2^I$ the probability that
\begin{displaymath} (1-\varepsilon)^{|X|} - \delta \leq
  \frac{\left|\bigcap_{s \in X} (C+s)\right|}{2^{|I|}} 
\leq (1-\varepsilon)^{|X|}+\delta \end{displaymath}
is at least 
$$1-
  2|G(X)|e^{-2^{|I|-2} \delta^2/|G(X)|} \geq 1-2^{|X|+1} e^{-2^{|I|-|X|-2}\delta^2}.$$
The probability that it happens for every $X$ of size $\leq m$ is at least
\begin{displaymath}1-2^{|I|\cdot (m+1)^2}\cdot 
  e^{- 2^{|I|-m-2}\delta^2} .\end{displaymath} 
If $m$ and $ \delta $ are fixed then this expression  approaches $1$ as
$|I|$ goes to infinity, 
since $\lim_{x \rightarrow \infty} P(x)e^{-x} = 0$ for any
  polynomial $P(x)$. 
It follows that for sufficiently large $|I|$ the probability that the
``random'' set $C$ has the required properties is $>0$. Thus
there exists an actual $C$ with these properties as well.~$\QED$

\section{Parameters of the construction}\label{four}
We will define now all the parameters
of the construction. The
actual relations (P1--P7 below) between these parameters make sense
only in the context
of the computations in which they are used, and 
are tailored to simplify the calculations in the following
sections.  The reason why we collected these definitions
 here is that there are many of them and 
the order in which
they are defined  is quite important. Nevertheless this section serves 
only as a reference.

The following notation will be used in the sequel. 
\begin{definition}
Suppose that $s : 
\omega \times \omega \longrightarrow \omega $.

Let $s^{(0)}(i,j)=i$ and $s^{(n+1)}(i,j)=s(s^{(n)}(i,j),j)$. 
Given $N \in \omega+1$, $n \in \omega $ and $f \in \omega^N$ let 
$$s^{(n)}(f)=\left\{\lft2(i,s^{(n)}\lft1(f(i),i\rgt1)\rgt2): i <N\right\}.$$
We will write $ s(f)$ instead of $ s^{(1)}(f) $.
\end{definition}

We define real
sequences   
$\{\varepsilon_i, \delta_i, \epsilon_i: i \in \omega \}$, intervals $\{I_i: i
\in \omega \}$, sets $\{C_i: i\in \omega \}$ and integers $\{m_i: i
\in \omega \}$. In addition we will define  functions  
$\bar{\sr},\tilde{\sr},\sr : \omega  \times
\omega \longrightarrow \omega $. 
The sequence $\{\varepsilon_i: i \in \omega \}$ is defined first. We require that
\begin{enumerate}
\item[(P1)] $0<\varepsilon_{i+1}< \varepsilon_i$ for $ i \in \omega $,
\item[(P2)] $\sum_{i \in \omega } \varepsilon_i< 1/2$.
\end{enumerate}
Set $ \epsilon_0=\delta_0=1$, $I_0=C_0=\emptyset, m_0=0$ and
$\bar{\sr}(n,0)=\tilde{\sr}(n,0)= \sr(n,0)=0$ for 
all $n \in \omega $.
Suppose that $\{\delta_i,\epsilon_i, I_i, C_i,m_i: i<N\}$ are
defined. Also assume that $\bar{\sr}(n,i),\tilde{\sr}(n,i)$ and 
$ \sr(n,i)$ are defined for $i<N$ and $ n \in
\omega $.

Put  $v_N=\left|\prod_{k<N}
    2^{I_k}\right|$, $l_N=\prod_{k<N} v_k$ and
define $\epsilon_N$ such that that
\begin{enumerate}
\item[(P3)] $0<v_N \cdot \epsilon_N \leq
  \varepsilon_N$,
\item[(P4)] $2^{l_N +N+2} \cdot\epsilon_N <  \epsilon_{N-1}$.
\end{enumerate}
Given $ \varepsilon_N$ and $\epsilon_N$ we will define for $k \in
\omega $
$$\bar{\sr}(k,N)=\left\{
  \begin{array}{ll}
\max\left\{l: \dfrac{k}{l+1} \epsilon_N^2
  \varepsilon_N^{l} > 4\right\} & \text{if }
k\epsilon_N^2 > 4\\
0 & \text{otherwise}
  \end{array}\right. .$$ 

Next let $\tilde{\sr}(k,N)=\bar{\sr}^{(2u_N)}(k,N)$, where
$u_N $ is the smallest integer $\geq \log_2(8/\epsilon_N^2)$.
Finally  
define
$$\sr(k,N)=\tilde{\sr}^{(2v_N+1)}(k,N).$$
Note that the functions $\bar{\sr}(\cdot,N),\tilde{\sr}(\cdot,N)$, and
$\sr(\cdot,N)$ are nondecreasing and unbounded.

Define
\begin{enumerate}
\item[(P5)] $m_{N}=\min\left\{m: \sr^{(N\cdot
      l_N)}(m,N)>0\right\}$,
\item[(P6)] $ \delta_N=2^{-N-2}\cdot \varepsilon_N^{m_N}$.
\end{enumerate}
Finally use \ref{choosec} to define $I_N$ and $C_N \subseteq
2^{I_N}$ for  $ \delta=\delta_N$,
$\varepsilon=\varepsilon_N$ and $m=m_N$. 

In addition we require that
\begin{enumerate}
\item[(P7)] $I_i$ are pairwise disjoint.
\end{enumerate}

The set $H$ that will witness that $\SM$ is not an ideal is defined as
$$H=\{x \in 2^\omega : \exists^\infty k \ x \rest I_k \not\in C_k\}.$$
Note that
$$\mu(H)\leq \mu\left(\bigcap_{n} \bigcup_{k>n} \{x \in 2^\omega : x
  \rest I_k \not\in C_k\}\right) \leq \sum_{k>n} \varepsilon_k + \delta_k
\stackrel{n \rightarrow
    \infty}{\longrightarrow}  0.$$

\section{More combinatorics}\label{five}
This section contains the core of the proof of \ref{newtrick}. This is
Theorem \ref{crucialgeneral} which 
is in the realm of finite combinatorics and 
concerns properties of the counting measure on finite product spaces. 
We will use the following notation:
\begin{definition}
  Suppose that $N_0<N \leq\omega $. Define
$\FF^N$ to be the collection of all sequences $\bF=\<f_i: i<N\>$
 such that $f_i \in \FF^{I_i}$ for $i < N$.
For $\bF \in \FF^N$ and $h \in \omega^N$, let
$$\FF^N_{\bF, h}=\left\{\bG\in \FF^N: \forall i < N \ \bG(i) \in
\FF^{I_i}_{\bF(i),h(i)}\right\}.$$
Similarly,
$$\FF^{N_0,N}_{\bF, h}=\left\{\bG \in \FF^N_{\bF, h}: \bG \rest N_0 =
  \bF \rest N_0\right\}.$$

We always require that for all $i < N$,
$$ \left|\dom\lft1(\bF(i)\rgt1)\right|+h(i) \leq m_i.$$

Let $\bC=\<C_i:i < \omega \>$ be the sequence of sets defined earlier.
For $N_0 <N$ and $\bF \in \FF^N$ let
$$(\bC)^{\bF}_{N_0}=\prod_{N_0\leq i< N} (C_i)^{\bF(i)} = \left\{s
  \in 2^{I_{N_0} \cup \dots \cup I_{N-1}}: \forall i \in [N_0, N)
 \ s \rest I_i \in (C_i)^{\bF(i)}\right\}.$$
We will write $(\bC)^{\bF}$ instead of $(\bC)^{\bF}_0$ and
$(C_{N-1})^{\bF(N-1)}$ instead of $(\bC)^{\bF}_{N-1}$.
\end{definition}

\begin{definition}
Suppose that $X$ is a finite set. A distribution is a function $m: X
\longrightarrow \reals$ 
such that $$0 \leq m(x) \leq \frac{1}{|X|}.$$
Define $\alpha_m$ to be the largest number $\alpha $ such that
$m'=\alpha \cdot m$ is a distribution,
and put 
$\wh{m}=\sum_{x \in X} m(x)$ and $\overline{\overline{m}}=\alpha_m \cdot \wh{m}$.

Suppose that a distribution $m$ on $X$ is given and $Y \subseteq X$. Define
$m_Y: Y \longrightarrow \reals^+$ as 
$$m_Y(x) = \frac{|X|}{|Y|} \cdot m(x).$$
\end{definition}
Note that 
$$\alpha_m=\frac{1}{|X|\cdot \max\{m(x): x\in X\}}.$$
Observe also that $(m_Y)_Z=m_Z$ if $Z \subseteq Y \subseteq X$.

A prototypical example of a distribution is defined as follows. Suppose
that $p \subseteq 2^\omega $ is a closed (or just measurable) set and $n \in \omega $. Let
$m$ be defined on $2^n$ as
$$m(s)=\mu(p \cap [s]) \text{ for } s \in 2^n.$$
Note that $\wh{m}=\mu(p)$.

The following lemmas list some easy observations concerning these
notions.
\begin{lemma}\label{lem1}
  Suppose that $N \in \omega $, $k^0+k_0 \leq m_N$, $f \in \FF^{I_N}_{\emptyset, k^0}$ and 
  $m$ is a distribution on $2^{I_N}$. There exist $f_0, f_1 \in
  \FF^{I_N}_{f, k_0}$ such that  $|f_0 \setminus f|=|f_1\setminus
  f|=k_0$ and
   $$\wh{m_{(C_N)^{f_0}}} \leq \wh{m_{(C_N)^{f}}} \leq \wh{m_{(C_N)^{f_1}}}.$$
\end{lemma}
\Proof 
For each $x \in 2^{I_N}$ and $h \in \FF^{I_N}_{\emptyset, k}$, let 
$h^0_x = h \cup \{(x,0)\}$ and $h^1_x = h \cup \{(x,1)\}$. Note that
there is $i \in \{0,1\}$ such that 
$$\wh{m_{(C_N)^{h^i_x}}} \leq \wh{m_{(C_N)^{h}}} \leq
\wh{m_{(C_N)^{h^{1-i}_x}}}.$$
Iteration of this procedure $k_0$ times will produce the required
examples.~$\QED$

\begin{lemma}\label{easy1}
  Suppose that $N_0 \leq N$ are natural numbers, $h^0,h_0 \in
\prod_{i<N} m_i$ satisfy $h_0(i)+h^0(i) \leq m_i$ for $i <N$, $\bF \in
 \FF^N_{\emptyset, h^0}$ and 
 $m$ is a distribution on
$2^{I_0 \cup \dots \cup I_{N-1}}$. 
Suppose that for every $\bG \in \FF^{N_0,N}_{\bF, h_0}$,
$a \leq \wh{m_{(\bC)^{\bG}}} \leq b$. Let $\bG^\star \in
\FF^{N_0,N}_{\bF, h_0}$ be such that
$|\dom\lft1(\bG^\star(i)\rgt1) \setminus \dom(\bF(i))|=h_1(i)<h_0(i)$
for $N_0 \leq i < N$. Then 
$$\forall \bG \in \FF^{N_0,N}_{\bG^\star, h_0-h_1}\
a \leq \wh{m_{(\bC)^{\bG}}} \leq b.$$
\end{lemma}
\Proof
Since $\FF^{N_0,N}_{\bG^\star, h_0-h_1} \subseteq
\FF^{N_0,N}_{\bF, h_0}$, the lemma is obvious.~$\QED$ 

The following theorem is a good approximation of the combinatorial
result that we require for the proof of \ref{newtrick}. 
The proof of it will give us a slightly stronger
but more technical result \ref{crucialgeneralmore}, which is
precisely what we need.
\begin{theorem}\label{crucialgeneral}
Suppose that $N_0 < N$ are natural numbers, $h^0,h_0 \in
\prod_{i<N} m_i$ satisfy $h_0(i)+h^0(i) \leq m_i$ for $i <N$, $\bF \in
 \FF^N_{\emptyset, h^0}$ and 
 $m$ is a distribution on
$2^{I_0 \cup \dots \cup I_{N-1}}$ such that 
$$\wh{m_{(\bC)^{\bF}}} \geq \frac{2\cdot \sum_{i=N_0}^N
  \epsilon_i}{\prod_{i=N_0}^N (1-8\epsilon_i)}.$$
There exists $\bF^\star \in
\FF^{N_0,N}_{\bF, h_0-\sr(h_0)}$ such that 
$$ \forall \bG \in \FF^{N_0,N}_{\bF^\star, \sr(h_0)}\ 
\wh{m_{(\bC)^{\bG}}} \geq \wh{m_{(\bC)^{\bF}}}\cdot\prod_{i=N_0}^{N-1} (1-8\epsilon_i)^2 -
\sum_{i=N_0}^{N-2} \epsilon_i.$$ 
\end{theorem}

{\bf Remark}.
It is worth noticing that the complicated formulas appearing in the
statement of this theorem are chosen to simplify the inductive  proof. Putting
them  aside,
the theorem can be formulated as follows: if $\wh{m_{(\bC)^{\bF}}}$ is 
sufficiently big (where big means only slightly larger than zero), then
there exists $\bF^\star \in
\FF^{N_0,N}_{\bF, h_0-\sr(h_0)}$ such that 
 for all $\bG \in \FF^{N_0,N}_{\bF^\star, \sr(h_0)}$ the value of 
$\dfrac{\wh{m_{(\bC)^{\bG}}}}{\wh{m_{(\bC)^{\bF}}}}$ cannot be significantly
smaller than~$1$.

The proof of \ref{crucialgeneral} 
will proceed   by induction on $N \geq N_0$, and the
following theorem 
corresponds to  the single
induction step.

Suppose that  $N \in \omega $ is fixed.

\begin{theorem}\label{comb}
If $k^0+k_0 \leq m_N$, $m$ is a distribution on $2^{I_N}$ and
$f \in \FF^{I_N}_{\emptyset, k^0}$ is such that $\overline{\overline{{m}_{({C_N})^{f}}}} \geq
2\epsilon_N$ then
there exists $f^\star \in
\FF^{I_N}_{f, k_0-\tilde{\sr}(k_0,N)}$ such that 
$$ \forall g \in \FF^{I_N}_{f^\star, \tilde{\sr}(k_0,N)}\ 
\wh{{m}_{({C_N})^{f^\star}}}\cdot (1+2\epsilon_N) \geq \wh{{m}_{({C_N})^{g}}} \geq \wh{{m}_{({C_N})^{f^\star}}}\cdot
(1- 2\epsilon_N),$$
and
$$ \forall g \in \FF^{I_N}_{f^\star, \tilde{\sr}(k_0,N)}\ 
\wh{{m}_{({C_N})^{g}}} \geq \wh{{m}_{({C_N})^{f}}}\cdot 
(1- 2\epsilon_N).$$

\end{theorem}
\Proof
We start with the following observation:
\begin{lemma}\label{small}
Suppose that $\wh{{m}_{({C_N})^f}} \geq \epsilon_N$.
  There exists $\tilde{f} \in
\FF^{I_N}_{f, k_0-\bar{\sr}(k_0,N)}$ such that 
$$ \forall g \in \FF^{I_N}_{\tilde{f}, \bar{\sr}(k_0,N)}\ 
\wh{{m}_{({C_N})^g}} \geq \wh{{m}_{({C_N})^f}}\cdot (1- \epsilon_N).$$
Similarly,
there exists $\tilde{f} \in
\FF^{I_N}_{f, k_0-\bar{\sr}(k_0,N)}$ such that 
$$ \forall g \in \FF^{I_N}_{\tilde{f}, \bar{\sr}(k_0,N)}\ 
\wh{{m}_{({C_N})^g}} \leq \wh{{m}_{({C_N})^{f}}}\cdot(1+ \epsilon_N).$$
\end{lemma}
\Proof
We will show only the first part, the second part is proved in the
same way. 
If $\bar{\sr}(k_0,N)=0$ then the lemma follows readily from
\ref{lem1}. 
Thus, suppose that $\bar{\sr}(k_0,N)>0$ and let 
$m_{({C_N})^f}$ be a distribution satisfying the requirements of the
lemma.

Construct,
by induction, a sequence $\{f_n : n < n^\star \}$
such that 
\begin{enumerate}
\item $f_0=f$,
\item $f_{n+1} \in \FF^{I_N}_{f_n,\bar{\sr}(k_0,N)}$,
\item $\wh{{m}_{({C_N})^{f_n}}} \geq
  \wh{{m}_{({C_N})^{f}}}\cdot\left(1+
    \dfrac{n}{2}\epsilon_N\varepsilon_N^{\bar{\sr}(k_0,N)}\right).$ 
\end{enumerate}
First notice that
$\bar{\sr}(k_0,N)$ was defined in such a way that 
\begin{multline*}\wh{{m}_{({C_N})^{f}}}\cdot\left(1+
    \frac{1}{2}\left(\frac{k_0}{\bar{\sr}(k_0,N)}-2\right)
\epsilon_N\varepsilon_N^{\bar{\sr}(k_0,N)}\right) \geq  \\   
\epsilon_N\left(1+
    \frac{1}{2}\left(\frac{k_0}{\bar{\sr}(k_0,N)}-2\right)
\epsilon_N\varepsilon_N^{\bar{\sr}(k_0,N)}\right)
\geq \\
\frac{1}{2}\frac{k_0}{\bar{\sr}(k_0,N)}\epsilon_N^2\varepsilon_N^{\bar{\sr}(k_0,N)} 
>1.
\end{multline*}
Therefore, after fewer than  $\dfrac{k_0}{\bar{\sr}(k_0,N)}-2 $ steps
the construction has to terminate  (otherwise $\wh{{m}_{({C_N})^{g}}}
>1 $ for some $g$, which is impossible).

Suppose that $f_n$ has been constructed. 

{\sc Case 1}. $\forall h \in \FF^{I_N}_{f_n, \bar{\sr}(k_0,N)} \
\wh{{m}_{({C_N})^{h}}} \geq \wh{{m}_{({C_N})^{f}}} \cdot(1-\epsilon_N)$. 
In this case put $\tilde{f}=f_n$ and finish the construction.
Observe that
\begin{multline*}
|\tilde{f}|+\bar{\sr}(k_0,N) \leq k^0+   n^\star\cdot
\bar{\sr}(k_0,N)+\bar{\sr}(k_0,N) \leq \\
k^0 +
\left(\frac{k_0}{\bar{\sr}(k_0,N)}-2\right) \cdot \bar{\sr}(k_0,N) +\bar{\sr}(k_0,N)
\leq k^0+k_0 - \bar{\sr}(k_0,N) < m_N.
\end{multline*}

{\sc Case 2}. $ \exists  h \in \FF^{I_N}_{f_n, \bar{\sr}(k_0,N)} \
\wh{{m}_{({C_N})^{h}}} < \wh{{m}_{({C_N})^{f}}} \cdot(1-\epsilon_N)$.
Using \ref{lem1} we can assume that $|h|=|f_n|+\bar{\sr}(k_0,N)$.

Consider the partition of $({C_N})^{f_n}$ given by $h$, i.e.
$$(C_N)^{f_n}=({C_N})^{h} \cup \lft1(({C_N})^{f_n} \setminus
({C_N})^{h}\rgt1).$$
Note that by considering the worst case we get 
\begin{multline*}
  \frac{|({C_N})^{h}|}{\left|({C_N})^{f_n}
      \setminus ({C_N})^{h}\right|} \geq \\
  \frac{(1-\varepsilon_N)^{m^0_{f_n}}\varepsilon_N^{m^1_{f_n}} \cdot
    \varepsilon_N^{\bar{\sr}(k_0,N)} - \delta_N
      }{(1-\varepsilon_N)^{m^0_{f_n}}\varepsilon_N^{m^1_{f_n}}+\delta_N -
      \left((1-\varepsilon_N)^{m^0_{f_n}}\varepsilon_N^{m^1_{f_n}}\cdot
\varepsilon_N^{\bar{\sr}(k_0,N)} - 
      \delta_N\right)} \geq \\
 \frac{\varepsilon_N^{\bar{\sr}(k_0,N)} -
    \dfrac{\delta_N}{(1-\varepsilon_N)^{m^0_{f_n}}
      \varepsilon_N^{m^1_{f_n}}}
    }{1-\varepsilon_N^{\bar{\sr}(k_0,N)}+\dfrac{2\delta_N}{(1-\varepsilon_N)^{m^0_{f_n}}
      \varepsilon_N^{m^1_{f_n}}}}.
\end{multline*}
Moreover, since $\delta_N \leq \dfrac{1}{2} \varepsilon_N^{m_N}$,  we have
$$\varepsilon_N^{\bar{\sr}(k_0,N)} \geq
2\cdot\frac{\delta_N}{(1-\varepsilon_N)^{m^0_{f_n}}\cdot
  \varepsilon_N^{m^1_{f_n}}},$$
and thus
$$\frac{\varepsilon_N^{\bar{\sr}(k_0,N)} -
    \dfrac{\delta_N}{(1-\varepsilon_N)^{m^0_{f_n}}
      \varepsilon_N^{m^1_{f_n}}}
    }{1-\varepsilon_N^{\bar{\sr}(k_0,N)}+\dfrac{2\delta_N}{(1-\varepsilon_N)^{m^0_{f_n}}
      \varepsilon_N^{m^1_{f_n}}}} \geq \dfrac{1}{2} \varepsilon_N^{\bar{\sr}(k_0,N)}.$$

It follows that
\begin{multline*} 
  \frac{1}{\wh{{m}_{({C_N})^{f}}}}\cdot \wh{{m}_{({C_N})^{f_n} \setminus ({C_N})^{h}}}
\geq  \left(1+\frac{n}{2} \epsilon_N
 \varepsilon_N^{\bar{\sr}(k_0,N)}\right)\cdot\frac{|({C_N})^{f_n}|}{\left|({C_N})^{f_n}
     \setminus ({C_N})^{h}\right|} -\\
(1-\epsilon_N)\cdot\frac{|({C_N})^{h}|}{\left|({C_N})^{f_n}
\setminus ({C_N})^{h}\right|}=
\end{multline*}
\begin{multline*}
\frac{|({C_N})^{f_n}|}{\left|({C_N})^{f_n}
\setminus ({C_N})^{h}\right|}+ \frac{n}{2} \epsilon_N \varepsilon_N^{\bar{\sr}(k_0,N)}
\cdot \frac{|({C_N})^{f_n}|}{\left|({C_N})^{f_n}
\setminus ({C_N})^{h}\right|}- \frac{|({C_N})^{h}|}{\left|({C_N})^{f_n}
\setminus ({C_N})^{h}\right|}+ \\
\epsilon_N  \frac{|({C_N})^{h}|}{\left|({C_N})^{f_n}
\setminus ({C_N})^{h}\right|}=
\end{multline*}
\begin{multline*}
  1+\frac{n}{2} \epsilon_N \varepsilon_N^{\bar{\sr}(k_0,N)} \cdot
  \frac{|({C_N})^{f_n}|}{\left|({C_N})^{f_n} \setminus
      ({C_N})^{h}\right|} + \epsilon_N
  \frac{|({C_N})^{h}|}{\left|({C_N})^{f_n}
      \setminus ({C_N})^{h}\right|} \geq\\
  1+ \frac{n}{2}\epsilon_N\varepsilon_N^{\bar{\sr}(k_0,N)} +
  \epsilon_N \cdot \frac{\varepsilon_N^{\bar{\sr}(k_0,N)} -
    \dfrac{\delta_N}{(1-\varepsilon_N)^{m^0_{f_n}}
      \varepsilon_N^{m^1_{f_n}}}
    }{1-\varepsilon_N^{\bar{\sr}(k_0,N)}+\dfrac{2\delta_N}{(1-\varepsilon_N)^{m^0_{f_n}}
      \varepsilon_N^{m^1_{f_n}}}}
  \geq \\
  1+ \frac{n}{2}\epsilon_N\varepsilon_N^{\bar{\sr}(k_0,N)}+
  \frac{1}{2}\epsilon_N\varepsilon_N^{\bar{\sr}(k_0,N)} \geq   
  1+\frac{n+1}{2} \epsilon_N \varepsilon_N^{\bar{\sr}(k_0,N)}.
\end{multline*}

Let $\{h_1, \dots, h_{2^{\bar{\sr}(k_0,N)}}\}$ be the list of all functions in
$\FF^{I_N}$ such that $\dom(h_i)=\dom(h) \setminus \dom(f_n)$.
Without loss of generality we can assume $f_n \cup h_1=h$.
The sets $({C_N})^{f_n \cup h_2}, \dots, ({C_N})^{f_n \cup h_{2^{\bar{\sr}(k_0,N)}}}$
  define a partition of the set $({C_N})^{f_n} \setminus ({C_N})^{h}$.
Since
$$\wh{{m}_{({C_N})^{f_n} \setminus ({C_N})^{h}}}
\geq \wh{{m}_{({C_N})^{f}}}\cdot\left(1+\frac{n+1}{2}\epsilon_N \varepsilon_N
^{\bar{\sr}(k_0,N)}\right)$$
it follows that there exists $2 \leq \ell \leq 2^{\bar{\sr}(k_0,N)}$ such that 
$$\wh{{m}_{({C_N})^{f_n \cup h_{\ell}}}} \geq
 \wh{{m}_{({C_N})^{f}}}\cdot\left(1+\frac{n+1}{2}\epsilon_N \varepsilon_N
^{\bar{\sr}(k_0,N)}\right).$$
Let $f_{n+1}=f_n \cup h_{\ell}$. This completes the induction.~$\QED$

{\sc Proof of \ref{comb}}. 
Suppose that $\overline{\overline{{m}_{({C_N})^{f}}}} =a_0\geq 2\epsilon_N$. Without
loss of generality we can assume that 
$\alpha_{m_{({C_N})^f}}=1$, that is
$\overline{\overline{{m}_{({C_N})^f}}}=\wh{{m}_{({C_N})^f}}$. This is because if we
succeed in proving the theorem for the distribution
$\alpha_{m_{({C_N})^f}} \cdot {m}_{({C_N})^f}$ then it must be true for
${m}_{({C_N})^f}$ as well.

Apply \ref{small} to get $f'
\in \FF^{I_N}_{f,k_0-\bar{\sr}(k_0,N)}$ such that 
$$ \forall g \in \FF^{I_N}_{f',\bar{\sr}(k_0,N)}\ 
\wh{{m}_{({C_N})^{g}}} \geq a_0(1-\epsilon_N).$$
Let 
$u_N$ be the smallest integer greater than $\log_2(8/\epsilon_N^2)$ and 
define by
induction sequences $\{f_i,a_i,b_i: i \leq u_N\}$ such that 
\begin{enumerate}
\item $b_0=1$ and $f_0=f'$,
\item $a_i,b_i \in \reals$ for $i \leq u_N$,
\item $|b_i-a_i| \leq 2^{-i}$ for $i \leq u_N$,
\item $f_{i+1} \in
  \FF^{I_N}_{f_i,\bar{\sr}^{(i+1)}(k_0,N)-\bar{\sr}^{(i+2)}(k_0,N)}$  for $i<u_N$,
\item $ \forall g \in \FF^{I_N}_{f_i,\bar{\sr}^{(i+1)}(k_0,N)}\
  a_i(1-\epsilon_N) \leq \wh{{m}_{({C_N})^{g}}} \leq b_i(1+\epsilon_N).$ 
\end{enumerate}

Suppose that $a_i,b_i$ and $f_i$ are defined and let 
$c = \wh{m_{(C_N)^{f_i}}}$. Observe that $c \geq
a_0\cdot(1-\epsilon_N) > \epsilon_N$.

If $|c-a_i| \leq 2^{-i-1}$ then
let $a_{i+1}=a_i$ and $b_{i+1}=c$. Apply \ref{small} to get $f_{i+1}
\in \FF^{I_N}_{f_i,\bar{\sr}^{(i+1)}(k_0,N)-\bar{\sr}^{(i+2)}(k_0,N)}$ such that 
$$ \forall g \in \FF^{I_N}_{f_{i+1},\bar{\sr}^{(i+2)}(k_0,N)}\ 
\wh{{m}_{({C_N})^{g}}} \leq b_{i+1}(1+\epsilon_N).$$
Otherwise let $a_{i+1}=c$ and $b_{i+1}=b_i$ and let $f_{i+1}
\in \FF^{I_N}_{h,\bar{\sr}^{(i+1)}(k_0,N)-\bar{\sr}^{(i+2)}(k_0,N)}$ be such that 
$$ \forall g \in \FF^{I_N}_{f_{i+1},\bar{\sr}^{(i+2)}(k_0,N)}\ 
\wh{{m}_{({C_N})^{g}}} \geq a_{i+1}(1-\epsilon_N).$$
Put $f^\star = f_{u_N}$. Note that
by the choice of $u_N$, $|b_{u_N} -a_{u_N}| \leq
\epsilon_N^2/8$. In addition,
$\bar{\sr}^{(u_N+1)}(k_0,N)> \bar{\sr}^{(2u_N)}(k_0,N)=\tilde{\sr}(k_0,N)$. 
Since $\wh{m_{(C_N)^{f^\star}}}$ is equal to either $a_{u_N}$ or
$b_{u_N}$, and $a_{u_N} \geq \varepsilon_N$, a simple computation
shows that for every $g \in \FF^{I_N}_{f^\star,\tilde{\sr}(k_0,N)}$, 
$$\wh{m_{(C_N)^{f^\star}}}\cdot(1-2\epsilon_N) \leq a_{u_N}(1-\epsilon_N) \leq
\wh{{m}_{({C_N})^{g}}} \leq b_{u_N}(1+\epsilon_N) \leq
\wh{m_{(C_N)^{f^\star}}}\cdot(1+2\epsilon_N),$$
and
$$\wh{{m}_{({C_N})^{g}}} \geq a_{u_N}(1-\epsilon_N) \geq
a_{0}(1-2\epsilon_N)=\wh{m_{(C_N)^{f}}}\cdot(1-2\epsilon_N).~\QED$$  

\bigskip

Before we start proving  \ref{crucialgeneral} we need to prove several
facts concerning distributions. The following notation will be used in 
the sequel.

\begin{enumerate}
\item $v_k = \left|2^{I_0 \cup \dots
  \cup I_{k-1}}\right|$ for $k \in \omega $. 
\item If $\bF \in \FF^{N}$ and $k<N$ then let $w_k(\bF)=
  \left|(\bC)^{\bF \rest k}\right|$.
\end{enumerate}

Suppose that $\bF \in \FF^{N+1}$ and $m$ is a distribution on $2^{I_0 \cup \dots
  \cup I_{N}}$. 
\begin{enumerate}
\item Let  $m^+_{(\bC)^{\bF}}$  be the distribution on $(C_N)^{\bF(N)}$ given by 
$$m_{(\bC)^{\bF}}^+(s)=\sum \left\{m_{(\bC)^{\bF}}(t): s \subseteq t
  \in (\bC)^{\bF}
\right\} \text{ for } s \in (C_N)^{\bF(N)} .$$
\item For $N_0\leq N$ and $t \in (\bC)^{\bF\rest N_0}$, let
  $m_{(\bC)^{\bF}}^t$ be a distribution
on $(\bC)^{\bF}_{N_0}$ defined as 
$$m_{(\bC)^{\bF}}^t(s)=m_{(\bC)^{\bF}}(t^\frown s) \text{ for } s \in
(\bC)^{\bF}_{N_0}.$$
\item Let $m_{(\bC)^{\bF}}^-$ be the distribution on $(\bC)^{\bF\rest N}$ defined as
$$m_{(\bC)^{\bF}}^-(t)=\wh{m_{(\bC)^{\bF}}^t} \text{ for } t \in (\bC)^{\bF\rest N}.$$
\end{enumerate}




\begin{lemma}\label{new0}
Suppose that $N_0 \leq N$, $\bF\in \FF^{N+1}$ and  $\bG \in
\FF^{N_0,N+1}_{\bF,h}$ for some $h \in \omega^\omega $.
Then
$$\left(m_{(\bC)^{\bF}}^t\right)_{(\bC)^{\bG}_{N_0}}=m_{(\bC)^{\bG}}^t.$$
\end{lemma}
\Proof
Fix $t \in (\bC)^{\bF\rest N_0}=(\bC)^{\bG\rest N_0}$ and
observe that for $s \in (\bC)^{\bF}_{N_0}$,
\begin{multline*}
  \left(m_{(\bC)^{\bF}}^t\right)_{(\bC)^{\bF}_{N_0}}(s)=
\frac{|(\bC)^{\bF}_{N_0}|}{|(\bC)^{\bG}_{N_0}|} 
\cdot  m_{(\bC)^{\bF}}(t^\frown s)=\\
\frac{|(\bC)^{\bF}_{N_0}|}{|(\bC)^{\bG}_{N_0}|}  \cdot  
\frac{v_{N+1}}{w_{N+1}(\bF)} m(t^\frown s)=\frac{v_{N+1}}{w_{N+1}(\bG)} 
\cdot m(t^\frown s)=   
m_{(\bC)^{\bG}}^t(s).~\QED
\end{multline*}

\begin{lemma}\label{new0a}
Suppose that $\bF \in \FF^{N+1}$ and $\bG \in \FF^{N,N+1}_{\bF, h}$
for some $h \in \omega^\omega $.
Then
$$(m_{(\bC)^{\bF}}^+)_{(C_N)^{\bG(N)}}=m_{(\bC)^{\bG}}^+.$$
\end{lemma}
\Proof
Similar to the proof of \ref{new0}.~$\QED$

\begin{lemma}\label{new1}
Suppose that $N_0 \leq N$, $\bF\in \FF^{N+1}$ and $t \in
(\bC)^{\bF\rest N_0}$. Then
  $$\overline{\overline{m_{(\bC)^{\bF}}^t}} \geq w_{N_0}(\bF)\cdot \wh{m_{(\bC)^{\bF}}^t}.$$
\end{lemma}
\Proof
Note that 
$$w_{N_0}(\bF) \cdot m_{(\bC)^{\bF}}(t^\frown s) \leq 
\frac{w_{N_0}(\bF)}{w_{N+1}(\bF)} = \frac{1}{|(\bC)^{\bF}_{N_0}|} .~\QED$$

The next two lemmas will be crucial in the recursive computations of
distributions.

\begin{lemma}\label{new2}
  Suppose that $N_0 \leq N$, $\bF,\bG \in \FF^{N+1}$, $\bF\rest
  [N_0,N]=\bG \rest [N_0,N] $ and $t \in
  (\bC)^{\bF\rest N_0} \cap (\bC)^{\bG\rest N_0}$. Then
$$\alpha_{m_{(\bC)^{\bF}}^t}\cdot m_{(\bC)^{\bF}}^t =
\alpha_{m_{(\bC)^{\bG}}^t }\cdot m_{(\bC)^{\bG}}^t.$$  
In particular, if $\bF^\star \in \FF^{N_0,N+1}_{\bF,h}$ for some $h\in 
\omega^\omega $ then
$$\frac{\wh{m^t_{(\bC)^{\bF\rest N_0{}^\frown
        \bF^\star \rest [N_0,N] }}}}{\wh{m^t_{(\bC)^{\bF}}}}=\frac{\wh{m^t_{(\bC)^{\bG\rest
        N_0{}^\frown \bF^\star \rest [N_0,N]}}}}{\wh{m^t_{(\bC)^{\bG}}}}.$$
\end{lemma}
\Proof
Note that under the assumptions the distributions $m_{(\bC)^{\bF}}^t$
and $m_{(\bC)^{\bG}}^t$ have the same domain and 
the fraction $\dfrac{m(t^\frown s)}{m^t_{(\bC)^{\bF}}(s)}$ has the constant value
for both $\bF$ and $\bG$.~$\QED$

\begin{lemma}\label{new3}
  Suppose that $\bF \in \FF^{N+1}$ and $\bG \in \FF^{N+1}_{\bF,h}$ for 
  some $h \in \omega^\omega$.
Then
$$\wh{m_{(\bC)^{\bG}}}=\sum_{t \in (\bC)^{\bG\rest N}}
m^-_{(\bC)^{\bG\rest N{}^\frown \bF(N)}}(t)\cdot
\frac{\wh{m^t_{(\bC)^{\bF\rest N{}^\frown \bG(N)}}}}{\wh{m^t_{(\bC)^{\bF}}}}.$$
\end{lemma}
\Proof
For $t \in (\bC)^{\bG\rest N}$,
$$\frac{m^-_{(\bC)^{\bG\rest N{}^\frown
      \bF(N)}}(t)}{\wh{m^t_{(\bC)^{\bF}}}}=
\frac{\displaystyle
\sum_{s' \in (C_N)^{\bF(N)}} \dfrac{v_{N+1}}{w_{N+1}\lft1(\bG\rest
     N{}^\frown \bF(N)\rgt1)} \cdot m(t^\frown s')
}
{\displaystyle\sum_{s' \in
   (C_N)^{\bF(N)}} \dfrac{v_{N+1}}{{w_{N+1}(\bF)}} \cdot m(t^\frown
    s')
}
=
\frac{w_N(\bF)}{w_N(\bG)} .$$
Therefore
\begin{multline*}
  \sum_{t \in (\bC)^{\bG\rest N}}
m^-_{(\bC)^{\bG\rest N{}^\frown \bF(N)}}(t)\cdot
\frac{\wh{m^t_{(\bC)^{\bF\rest N{}^\frown
        \bG(N)}}}}{\wh{m^t_{(\bC)^{\bF}}}} =
\sum_{t \in (\bC)^{\bG\rest N}} \frac{w_N(\bF)}{w_N(\bG)}\cdot
\wh{m^t_{(\bC)^{\bF\rest N{}^\frown 
        \bG(N)}}} =\\
\sum_{t \in (\bC)^{\bG\rest N}} \frac{w_N(\bF)}{w_N(\bG)}\cdot
\sum_{t \subseteq s \in (\bC)^{\bG}} \frac{v_{N+1}}{w_{N+1}\lft1(\bF\rest
  N{}^\frown \bG(N)\rgt1)}\cdot m(s) = \\
\sum_{t \in (\bC)^{\bG\rest
    N}} \sum_{t \subseteq s \in (\bC)^{\bG}}
\frac{v_{N+1}}{w_{N+1}(\bG)}\cdot m(s)=
\sum_{t \in (\bC)^{\bG\rest
    N}} \sum_{t \subseteq s \in (\bC)^{\bG}}
m_{(\bC)^{\bG}}(s)=
\wh{m_{(\bC)^{\bG}}}.\QED
\end{multline*}

We will need one more definition:
\begin{definition}
  Suppose that $m$ is a distribution on $X$ and $U \subseteq X$. 
Let $ \MU{U}$ be the distribution on $X$ defined as 
$$ \MU{U}(x)=\left\{
  \begin{array}{ll}
m(x) & \text{if }x \in U\\
0 &    \text{otherwise}
  \end{array}\right. \text{ for } x \in X.$$

\end{definition}

Now we are ready to prove theorem \ref{crucialgeneral}. 
For technical
reasons we will need  a somewhat stronger result stated below.

\begin{theorem}\label{crucialgeneralmore}
Suppose that $N_0 < N$ are natural numbers, $h^0,h_0 \in
\prod_{i<N} m_i$ satisfy $h_0(i)+h^0(i) \leq m_i$ for $i <N$, $\bF \in
 \FF^N_{\emptyset, h^0}$ and 
 $m$ is a distribution on
$2^{I_0 \cup \dots \cup I_{N-1}}$ such that 
$$\wh{m_{(\bC)^{\bF}}} \geq \frac{2 \sum_{i=N_0}^N
  \epsilon_i}{\prod_{i=N_0}^N (1-8\epsilon_i)}.$$
There exist $\bF^\star \in
\FF^{N_0,N}_{\bF, h_0-\sr(h_0)}$ and $U^\star \subseteq 2^{I_0 \cup \dots
  \cup I_{N-1}}$ such that 
$$\wh{ \left(\MU{U^\star}\right)_{(\bC)^{\bF^\star}}} \geq
\wh{m_{(\bC)^{\bF}}}\cdot \prod_{i=N_0}^{N-1} (1-8\epsilon_i) -
\sum_{i=N_0}^{N-2} \epsilon_i,$$
and 
for any  $\bG \in \FF^{N_0,N}_{\bF^\star, \sr(h_0)}$ and $t \in
(\bC)^{\bG\rest M_0}$, $M_0 \in [N_0,N)$,
$$\wh{\left(\MU{U^\star}\right)_{(\bC)^{\bG}}^t} \geq 
\wh{\left(\MU{U^\star}\right)_{(\bC)^{\bF^\star}}^t}\cdot\prod_{i=M_0}^{N-1}
(1-4\epsilon_i).$$   
\end{theorem}
\Proof First notice that \ref{crucialgeneral} follows from
\ref{crucialgeneralmore}.
If $\bF^\star$ and $U^\star$ are as required, then for all $\bG \in
\FF^{N_0,N}_{\bF^\star, \sr(h_0)}$, 
\begin{multline*}
\wh{m_{(\bC)^{\bG}}} \geq \wh{\left(\MU{U^\star}\right)_{(\bC)^{\bG}}} \geq
\sum_{t \in (\bC)^{\bG \rest N_0}}
\wh{\left(\MU{U^\star}\right)_{(\bC)^{\bG}}^t} \geq \\
\sum_{t \in (\bC)^{\bG \rest N_0}} \left(\wh{
\left(\MU{U^\star}\right)_{(\bC)^{\bF^\star}}^t}\cdot\prod_{i=N_0}^{N-1}
(1-4\epsilon_i)\right)\geq \\
\prod_{i=N_0}^{N-1}
(1-4\epsilon_i) \cdot \left(\sum_{t \in (\bC)^{\bF^\star \rest N_0}} 
\wh{\left(\MU{U^\star}\right)_{(\bC)^{\bF^\star}}^t}\right) =
 \prod_{i=N_0}^{N-1}
(1-4\epsilon_i) \cdot \wh{\left(\MU{U^\star}\right)_{(\bC)^{\bF^\star}}}\geq \\
\prod_{i=N_0}^{N-1}
(1-4\epsilon_i)\cdot 
\left(\wh{m_{(\bC)^{\bF}}}\cdot \prod_{i=N_0}^{N-1} (1-8\epsilon_i) -
\sum_{i=N_0}^{N-2} \epsilon_i \right) \geq \\
\wh{m_{(\bC)^{\bF}}}\cdot\prod_{i=N_0}^{N-1} (1-8\epsilon_i)^2 -
\sum_{i=N_0}^{N-2} \epsilon_i.
\end{multline*}

\bigskip

We will proceed by induction on $N$. 
If $N=N_0$ then the theorem  is trivially true.
Thus, suppose that the result holds for some $N \geq N_0$ and consider $N+1$.
Let 
$\bF \in \FF^{N_0,N+1}_{\emptyset, h^0}$ 
and let $m$ be a distribution on $2^{I_0 \cup \dots \cup I_{N}}$ such
that 
$$\wh{m_{(\bC)^{\bF}}}\geq \frac{2\sum_{i=N_0}^N
  \epsilon_i}{\prod_{i=N_0}^N(1-\epsilon_i)}.$$ 

Recall that by \ref{new0a}, 
$$\overline{\overline{m^+_{(\bC)^{\bF}}}}\geq \wh{m^+_{(\bC)^{\bF}}}=\wh{m_{(\bC)^{\bF}}}\geq 2\epsilon_N,$$
and apply \ref{comb} with
$m=m^+_{(\bC)^{\bF}}$, $k^0 = \left|\dom\lft1(\bF(N)\rgt1)\right|$,
$k_0=h_0(N)$ to get $\tilde{f}_0  
\in \FF^{I_N}_{\bF(N),k_0-\tilde{\sr}(k_0,N)}$ such that 
$$ \forall g \in \FF^{I_N}_{\tilde{f}_0,\tilde{\sr}(k_0,N)}\ 
\wh{{m^+}_{({\bC})^{\bF\rest N{}^\frown g}}} \geq
\wh{m_{(\bC)^{\bF}}}\cdot (1-2\epsilon_N).$$

Let $\{s_i: 1 \leq i \leq w_N(\bF)\}$ be an
enumeration of $(\bC)^{\bF\rest N}$. 
By induction, build  a sequence
$\{\tilde{f}_{i}: i \leq w_N(\bF)\}$ such that 
\begin{enumerate}
\item $\tilde{f}_i \subseteq \tilde{f}_{i+1}$,
\item $k_0 - |\dom(\tilde{f}_i)| \geq \tilde{\sr}^{(2i+1)}(k_0,N)$,
\item for every $i \geq 1$ one of the following conditions holds:
    \begin{enumerate}
    \item $ \forall g \in \FF^{I_N}_{\tilde{f}_i,\tilde{\sr}^{(2i+1)}(k_0,N) }\ 
\wh{{m^{s_i}}_{(\bC)^{\bF\rest N{}^\frown g}}} < \dfrac{2\epsilon_N}{w_N(\bF)}$,
\item for all $g \in \FF^{I_N}_{\tilde{f}_{i}, \tilde{\sr}^{(2i+1)}}$, 
$$\wh{{m^{s_i}}_{(\bC)^{\bF\rest N{}^\frown\tilde{f}_{i}}}} \cdot
(1- 2\epsilon_N) \leq \wh{{m^{s_i}}_{(\bC)^{\bF\rest N{}^\frown g}}} \leq
\wh{{m^{s_i}}_{(\bC)^{\bF\rest N{}^\frown\tilde{f}_{i}}}}  \cdot
(1+ 2\epsilon_N).$$ 
    \end{enumerate}
\end{enumerate}

\bigskip

Suppose that $\tilde{f}_i$ is given. 
If
$$ \forall g \in \FF^{I_N}_{\tilde{f}_i,\tilde{\sr}^{(2i+3)}(k_0,N) }\ 
\wh{{m^{s_{i+1}}}_{(\bC)^{\bF\rest N{}^\frown g}}} < \frac{2\epsilon_N}{w_N(\bF)}$$
then put $\tilde{f}_{i+1}=\tilde{f}_{i}$.

Otherwise, let $\tilde{f}_{i+1}' \in
\FF^{I_N}_{\tilde{f}_i,\tilde{\sr}^{(2i+3)}(k_0,N) }$ be chosen so
that  
$$\wh{{m^{s_{i+1}}}_{(\bC)^{\bF\rest N{}^\frown\tilde{f}_{i+1}'}}} \geq
\frac{2\epsilon_N}{w_N(\bF)}.$$
 In particular, by \ref{new1},
$\overline{\overline{{m^{s_{i+1}}}_{(\bC)^{\bF\rest N{}^\frown\tilde{f}_{i+1}'}}}}\geq 2\epsilon_N$.
Let $\tilde{k}=k_0 - |\dom{\tilde{f}_{i+1}'}|$.
By \ref{comb}, 
there exist $\tilde{f}_{i+1} \in
\FF^{I_N}_{\tilde{f}_{i+1}',\tilde{k} - \tilde{\sr}(\tilde{k},N)
  }$ such that for all $g \in \FF^{I_N}_{\tilde{f}_{i+1}, \tilde{\sr}(\tilde{k},N)}$,
$$ \wh{{m^{s_{i+1}}}_{(\bC)^{\bF\rest N{}^\frown\tilde{f}_{i+1}}}} \cdot
(1+ 2\epsilon_N) \geq \wh{{m^{s_{i+1}}}_{(\bC)^{\bF\rest N{}^\frown g}}} \geq
\wh{{m^{s_{i+1}}}_{(\bC)^{\bF\rest N{}^\frown\tilde{f}_{i+1}}}}  \cdot
(1- 2\epsilon_N).$$ 
Note that  $\tilde{k}\geq k_0 - 
|\dom{\tilde{f}_{i}}|-\tilde{\sr}^{(2i+3)}(k_0,N)$. Using the
induction hypothesis we get that
$\tilde{k} \geq \tilde{\sr}^{(2i+1)}(k_0,N)-\tilde{\sr}^{(2i+3)}(k_0,N)
\geq \tilde{\sr}^{(2i+2)}(k_0,N)$. It follows that
$\tilde{\sr}(\tilde{k},N) \geq \tilde{\sr}^{(2i+3)}(k_0,N)$ and 
${k}_0 - |\dom(\tilde{f}_{i+1})| \geq
\tilde{\sr}^{(2i+3)}(k_0,N)$, which finishes the induction.

\bigskip

 Let $\bF^\star(N) =\tilde{f}_{w_N(\bF)}$. Since $w_N(\bF) \leq \left|2^{I_0
     \cup \dots \cup I_{N-1}}\right|$ it follows that
 $\sr(k_0,N) \leq \tilde{\sr}^{(2w_N(\bF)+1)}(k_0,N)$. Thus
$\bF^\star(N) \in
\FF^{I_N}_{\bF(N),h_0(N)-\sr(k_0,N)}$.

Observe that $\wh{{m^{s}}_{(\bC)^{\bF\rest N{}^\frown
      g}}}=m^-_{(\bC)^{\bF\rest N{}^\frown g}}(s)$ for every $s \in
(\bC)^{\bF\rest N}$.
In particular, $\wh{{m^{s}}_{(\bC)^{\bF\rest N{}^\frown
      \bF^\star(N)}}}=m^-_{(\bC)^{\bF\rest N{}^\frown \bF^\star(N)}}(s)$.

By the construction, for every $g \in
\FF^{I_N}_{\bF^\star(N),\sr(k_0,N)}$ and $s \in (\bC)^{\bF\rest
  N}$,
$$ m^-_{(\bC)^{\bF\rest N{}^\frown
    \bF^\star(N)}}(s)\cdot \frac{1-2\epsilon_N}{1+2\epsilon_N}\leq
\wh{{m^{s}}_{(\bC)^{\bF\rest N{}^\frown g}}} \leq 
  m^-_{(\bC)^{\bF\rest N{}^\frown \bF^\star(N)}}(s)\cdot
  \frac{1+2\epsilon_N}{1-2\epsilon_N}  $$
or otherwise
$$\wh{{m^{s}}_{(\bC)^{\bF\rest N{}^\frown \bF^\star(N)}}} \leq
\frac{2\epsilon_N}{w_N(\bF)} \quad \text{ and }\quad
\wh{{m^{s}}_{(\bC)^{\bF\rest N{}^\frown g}}}\leq 
\frac{2\epsilon_N}{w_N(\bF)}.$$

Moreover, by the choice of $\tilde{f}_0$,  for every $g \in
\FF^{I_N}_{\bF^\star(N),\sr(k_0,N)}$,
$$\wh{{m^+}_{({\bC})^{\bF\rest N{}^\frown g}}} \geq
\wh{m_{(\bC)^{\bF}}}\cdot(1-2\epsilon_N).$$ 

Even though we do not have much control over the values of
$m^-_{(\bC)^{\bF\rest N{}^\frown \bF^\star(N)}}(s)$ we can show that
many  
of them are larger than $\dfrac{2\epsilon_N}{w_N(\bF)}$. 
Let
$$U=\left\{s \in (\bC)^{\bF\rest N}: m^-_{(\bC)^{\bF\rest
      N{}^\frown \bF^\star(N)}}(s) \geq
  \frac{2\epsilon_N}{w_N(\bF)}\right\}.$$ 
Note that for every $g \in
\FF^{I_N}_{\bF^\star(N),\sr(k_0,N)}$,
\begin{multline*}
  (1-2\epsilon_N) \cdot \wh{m_{(\bC)^{\bF}}} \leq
  \wh{{m^+}_{({\bC})^{\bF\rest N{}^\frown g}}} =  
\wh{{m}_{(\bC)^{\bF\rest N{}^\frown g}}} =\wh{{m^-}_{({\bC})^{\bF\rest
      N{}^\frown g}}} \leq\\
  \sum_{s \in U} 
m^-_{(\bC)^{\bF\rest N{}^\frown g}}(s) + \sum_{s \in
  ({\bC})^{\bF\rest N} \setminus U} m^-_{(\bC)^{\bF\rest N{}^\frown
    g}}(s) \leq \\ 
\frac{1+2\epsilon_N}{1-2\epsilon_N}\cdot\sum_{s \in U} m^-_{(\bC)^{\bF\rest N{}^\frown
    \bF^\star(N)}}(s) + w_N(\bF)\cdot \frac{2\epsilon_N}{w_N(\bF)}\cdot 
\frac{1+2\epsilon_N}{1-2\epsilon_N} \leq \\
2\epsilon_{N}\cdot \frac{1+2\epsilon_N}{1-2\epsilon_N} + \frac{1+2\epsilon_N}{1-2\epsilon_N}\cdot\sum_{s \in U} m^-_{(\bC)^{\bF\rest
    N{}^\frown \bF^\star(N)}}(s). 
\end{multline*}
It follows that
$$ \sum_{s \in U} m^-_{(\bC)^{\bF\rest N{}^\frown \bF^\star(N)}}(s)
\geq  \frac{(1-2\epsilon_N)^2}{1+2\epsilon_N}\cdot \wh{m_{(\bC)^{\bF}}} -
2\epsilon_N \geq (1-8\epsilon_N)\cdot \wh{m_{(\bC)^{\bF}}} -
2\epsilon_{N}.$$ 
Define distribution $m^\star$ on $2^{I_0\cup \dots \cup I_{N-1}}$ as 
$$m^\star(s) = \left\{
  \begin{array}{ll}
m^-_{(\bC)^{\bF\rest N{}^\frown \bF^\star(N)}}(s) & \text{if }
m^-_{(\bC)^{\bF\rest N{}^\frown \bF^\star(N)}}(s) \geq
\dfrac{2\epsilon_N}{w_N(\bF)}\\ 
0 & \text{otherwise}
  \end{array}\right. .$$

Clearly,
$$\wh{m^\star_{(\bC)^{\bF\rest N}}}=\sum_{s \in U} m^-_{(\bC)^{\bF\rest
    N{}^\frown \bF^\star(N)}}(s) \geq\\ 
 \wh{m_{(\bC)^{\bF}}}\cdot(1-8\epsilon_N) - 2\epsilon_{N}\geq 
\frac{2 \sum_{i=N_0}^{N-1}
  \epsilon_i}{\prod_{i=N_0}^{N-1} (1-8\epsilon_i)}.$$
Apply the induction hypothesis to $m^\star$, $\bF\rest N$ and $h_0\rest N$ to
obtain $\bF^\star\rest N$ and $V^\star$ as in \ref{crucialgeneralmore}.
Let 
$$U^\star = \left\{s \in 2^{I_0 \cup \dots \cup I_N}: s \rest I_0 \cup
  \dots \cup I_{N-1} \in V^\star 
\cap U\right\}.$$
It remains to check that $\bF^\star$ and $U^\star$ have the required
properties.
\begin{multline*}
\wh{\left(\MU{U^\star}\right)_{(\bC)^{\bF^\star}}} = \sum_{s\in (\bC)^{\bF^\star
    \rest N}} \left(\MU{U^\star}\right)^-_{(\bC)^{\bF^\star}} =\sum_{s\in (\bC)^{\bF^\star
    \rest N}}  \left(\MU{V^\star}^\star\right)_{(\bC)^{\bF^\star\rest N}}(s) =\\
\wh{ \left(\MU{V^\star}^\star\right)_{(\bC)^{\bF^\star\rest N}}} \geq
\wh{m^\star_{(\bC)^{\bF\rest N}}}\cdot \prod_{i=N_0}^{N-1} (1-8\epsilon_i) -
\sum_{i=N_0}^{N-2} \epsilon_i \geq\\
 \left(\wh{m_{(\bC)^{\bF}}}\cdot(1-8\epsilon_N) - 2\epsilon_{N}\right)
\cdot \prod_{i=M_0}^{N-1} (1-8\epsilon_i) - \sum_{i=M_0}^{N-2}
  \epsilon_i \geq\\
 \wh{m_{(\bC)^{\bF}}}\cdot \prod_{i=M_0}^{N}
  (1-8\epsilon_i) - \sum_{i=M_0}^{N-1} 
  \epsilon_i,
\end{multline*}
which gives the first condition.

To verify the second condition 
suppose that $\bG \in \FF^{N_0,N+1}_{\bF^\star, \sr(h_0)}$, $M_0
\in [N_0,N]$ and $t \in (\bC)^{\bG\rest M_0}$.
 By the inductive hypothesis we have that
$$\wh{\left(\MU{V^\star}^\star\right)^t_{(\bC)^{\bG\rest N}}} \geq 
\wh{\left(\MU{U^\star}^\star\right)_{(\bC)^{\bF^\star\rest
      N}}^t}\cdot\prod_{i=M_0}^{N-1} (1-4\epsilon_i).$$  

By \ref{new2} and \ref{new3},
\begin{multline*}
  \wh{\left(\MU{U^\star}\right)^t_{(\bC)^{\bG}}} = \sum_{t \subseteq s \in (\bC)^{\bG\rest N}}
\left(\MU{U^\star}\right)^-_{(\bC)^{\bG\rest N{}^\frown \bF^\star(N)}}(s)\cdot
\frac{\wh{\left(\MU{U^\star}\right)^s_{(\bC)^{\bF^\star\rest N{}^\frown
        \bG(N)}}}}{\wh{\left(\MU{U^\star}\right)^s_{(\bC)^{\bF^\star}}}}=\\
\sum_{t \subseteq s \in (\bC)^{\bG\rest N}}
\left(\MU{U^\star}\right)^-_{(\bC)^{\bG\rest N{}^\frown \bF^\star(N)}}(s)\cdot
\frac{\wh{\left(\MU{U^\star}\right)^s_{(\bC)^{\bF\rest N{}^\frown
        \bG(N)}}}}{\wh{\left(\MU{U^\star}\right)^s_{(\bC)^{\bF\rest N{}^\frown
        \bF^\star(N)}}}}=\\
\sum_{t \subseteq s \in (\bC)^{\bG\rest N}}
\left(\MU{U^\star}\right)^-_{(\bC)^{\bG\rest N{}^\frown \bF^\star(N)}}(t)\cdot
\frac{\wh{\left(\MU{U}\right)^s_{(\bC)^{\bF\rest N{}^\frown
        \bG(N)}}}}{\wh{\left(\MU{U}\right)^s_{(\bC)^{\bF\rest N{}^\frown
        \bF^\star(N)}}}}.
\end{multline*}
Now
\begin{multline*}
  \sum_{t \subseteq s \in (\bC)^{\bG\rest N}}
\left(\MU{U^\star}\right)^-_{(\bC)^{\bG\rest N{}^\frown \bF^\star(N)}}(s)\cdot
\frac{\wh{\left(\MU{U}\right)^s_{(\bC)^{\bF\rest N{}^\frown
        \bG(N)}}}}{\wh{\left(\MU{U}\right)^s_{(\bC)^{\bF\rest N{}^\frown
        \bF^\star(N)}}}} \geq \\
\sum_{t \subseteq s \in (\bC)^{\bG\rest N}}
\left(\MU{U^\star}^\star\right)^t_{(\bC)^{\bG\rest N}}(s)\cdot
\frac{1-2\epsilon_N}{1+2\epsilon_N} =\sum_{t \subseteq s \in (\bC)^{\bG\rest N}}
\left(\MU{V^\star}^\star\right)^t_{(\bC)^{\bG\rest N}}(s)\cdot
\frac{1-2\epsilon_N}{1+2\epsilon_N} \geq \\
\frac{1-2\epsilon_N}{1+2\epsilon_N}\cdot
\wh{\left(\MU{V^\star}^\star\right)_{(\bC)^{\bF^\star\rest 
      N}}^t}\cdot\prod_{i=M_0}^{N-1} (1-4\epsilon_i) 
 \geq \wh{\left(\MU{V^\star}^\star\right)_{(\bC)^{\bF^\star\rest
      N}}^t}\cdot\prod_{i=M_0}^{N} (1-4\epsilon_i) =\\
\wh{\left(\MU{U^\star}\right)_{(\bC)^{\bF^\star}}^t}\cdot\prod_{i=M_0}^{N} (1-4\epsilon_i),
\end{multline*}
which concludes the proof.

\section{Measures and norms}\label{six}
In this section we will examine the consequences of the combinatorial
results proved earlier on  measures on $2^\omega $.

For $U \subseteq  2^I$, $[U]=\{x \in
2^\omega : x \rest I \in  U\}$. 

If $p \subseteq 2^{<\omega}$
is a tree,  $s \in p$, and $N \in \omega $, then
\begin{enumerate}
\item $[p]$ denotes the set of  branches of $p$,
\item $p_s = \{t \in p: t \subseteq s \text{ or } s \subseteq t\}$,
\item $p^N= p \rest (I_0 \cup \dots \cup I_{N-1})$.
\end{enumerate}
We will identify product with concatenation, i.e., 
$(s,t) $ with  $s^\frown t$, and similarly for infinite
products.
Most of the time we will also identify $p$ with $[p]$. 

\begin{definition}
  Let $\mu_{(\bC)^{\Ff}}$ be the measure on $(\bC)^{\Ff}$ defined as the 
product of counting measures on the coordinates.
In other words, if $s \in 2^{I_k}$ then 
\begin{displaymath}\mu_{(\bC)^{\Ff}}([s])=\left\{
  \begin{array}{ll}
\left|(C_k)^{\Ff(k)}\right|^{-1} & \text{if } s \in (C_k)^{\Ff(k)}\\
0 & \text{otherwise}
  \end{array}\right. . \end{displaymath}
\end{definition}

Given a perfect set $p \in \perf$, 
$$\mu_{(\bC)^{\Ff}}(p)=\lim_{N \rightarrow \infty} \frac{\left|p^N \cap 
  (\bC)^{\Ff\rest N}\right|}{\left|(\bC)^{\Ff\rest N} \right|}.$$
Note that $\mu_{(\bC)^{\Ff}}(p)=\mu_{(\bC)^{\Ff}}\lft1(p\cap (\bC)^{\Ff}\rgt1)$.

\begin{definition}
For a function $f \in \omega^\omega $ define 
$\log_{\sr}(f) \in \omega^\omega $ as 
$$\log_{\sr}(f)(N)=\max\left\{k: \sr^{(k\cdot l_N)}(f(N),N)> 0\right\}.$$

For $h_1, h_s \in \omega^\omega $ define $h_1 \simeq h_2$ if
$\log_{\sr}(h_1)=\log_{\sr}(h_2)$. Clearly $\simeq $ is an equivalence 
relation. 

Let ${\mathcal X} $ be the collection of functions $f \in
\omega^\omega $ such that 
\begin{enumerate}
\item $\lim_{m \rightarrow \infty}\log_{\sr}(f)(m)=\infty,$
\item $f = \min\{g: f \simeq g\}$.
\end{enumerate}

For $f \in \omega^\omega $ define functions $\bar{f},f^- \in {\mathcal
  X} $ as
follows:
$\overline{f} = {\mathcal X} \cap \{g: f \simeq g\},$ and 
$$f^-(n)=\left\{
  \begin{array}{ll}
\min\{k: \log_{\sr}(k,n)=\log_{\sr}(f(n),n)-1\}& \text{if
  }\log_{\sr}(f(n),n)>0\\
0 & \text{otherwise}
  \end{array}\right.  .$$

If $f \in {\mathcal X}$ and $ n \in \omega $ let 
$i_f(n)=\max\{k: \log_{\sr}(f)(k) \leq n\}.$
\end{definition}
{\bf Remarks}.
Note that ${\mathcal X} \neq \emptyset$. By P5, $\overline{h} \in {\mathcal
  X} $, where 
$h(k)=m_k$ for $k\in \omega$.
Also, $\lim_{n \rightarrow \infty} i_f(n)=\infty$ for $f \in
{\mathcal X} $.
The purpose of the restriction put on the set ${\mathcal X} $ is to make 
the mapping $f \mapsto \log_{\sr}(f)$  one-to-one.
In practice, we will only use the fact that if $\log_{\sr}(f)(n)=0$
then $f(n)=0$.

\begin{definition}
  For a perfect set $p \subseteq 2^\omega $, $\Ff \in \FF^\omega$,
  $N\in \omega $ and
  $h \in {\mathcal X}  $, define
$$\bv{p, \Ff, h}_N=\inf \left\{\mu_{(\bC)^{\mathbf G}}(p) : {\mathbf G}
  \in \FF^{N,\omega}_{\Ff,h}\right\}.$$ 
We will write $\bv{p, \Ff, h}$ instead of $\bv{p, \Ff, h}_0$.
\end{definition}

The following easy lemma lists some basic properties of these notions.
\begin{lemma}\label{easy2}
  \begin{enumerate}
\item The sequence $\left\{\dfrac{\left|p^N \cap 
  (\bC)^{\Ff\rest N}\right|}{\left|(\bC)^{\Ff\rest N} \right|}: k 
\in \omega\right\}$ is  monotonically decreasing for every $p \in \perf$,
\item $\bv{p, \Ff_1,h_1}_N \geq \bv{p, \Ff_2, h_2}_N$ if $\Ff_1 \in
  \FF^{N,\omega}_{\Ff_2,h_2-h_1}$,
\item if $p_1 \cap p_2 = \emptyset$ then
$\bv{p_1\cup p_2, \Ff, h}_N \geq \bv{p_1, \Ff, h}_N +\bv{p_2, \Ff, h}_N$.
  \end{enumerate}
\end{lemma}
\Proof
(1) is obvious, and (2) follows from
\ref{easy1}.

(3) Take $ \varepsilon>0$ and let $\bG \in \FF^{N,\omega}_{\Ff,h}$ be
such that 
$$\bv{p_1\cup p_2, \Ff,
    h}_N+\varepsilon \geq \mu_{(\bC)^{\bG}}(p_1 \cup p_2).$$
Now
\begin{multline*}
\bv{p_1\cup p_2, \Ff,
    h}_N+\varepsilon \geq \mu_{(\bC)^{\bG}}(p_1 \cup p_2)\geq
    \mu_{(\bC)^{\bG}}(p_1 \cup p_2) \geq \\
  \mu_{(\bC)^{\bG}}(p_1)+\mu_{(\bC)^{\bG}}(p_2) \geq
 \bv{p_1, \Ff, h}_N
      +\bv{p_2, \Ff, h}_N.
\end{multline*}
Thus 
$\bv{p_1\cup p_2, \Ff,
    h}_N+\varepsilon \geq \bv{p_1, \Ff, h}_N
      +\bv{p_2, \Ff, h}_N$ and  the inequality follows.~$\QED$



The following two theorems are the key to the whole construction. 
\begin{theorem}\label{main}
  Suppose that  $\mu_{(\bC)^{\Ff}}(p)>0$, $h \in {\mathcal X}  $ and $0<\varepsilon <1$. 
Then there exist $p^\star \subseteq p$,  $h^\star \in {\mathcal X} $,
$N_0 \in \omega $ and $\Ff^\star \in 
\FF^{N_0, \omega}_{\bF, h-h^\star}$  such that 
$$\mu_{(\bC)^{\Ff^\star}}(p^\star) \geq
(1-\varepsilon)\cdot\mu_{(\bC)^{\Ff}}(p), \quad
\bv{p^\star,\bF^\star,h^\star} \geq
(1-2\varepsilon)\cdot\mu_{(\bC)^{\Ff}}(p)$$
and
$$\forall N \ \forall s \in (p^\star)^{N} \
\bv{p^\star_s,\bF^\star,h^\star}_N>0.$$
Moreover, we can require that $h^\star(N)=\overline{\sr(h)}(N)=h^-(N)$ for $N \geq N_0$.
\end{theorem}
\Proof
Find $N_0\in \omega $ such that 
\begin{enumerate}
\item $\mu_{(\bC)^{\bF}}(p)>\dfrac{2\sum_{i=N_0}^\infty
  \epsilon_i}{\prod_{i=N_0}^\infty (1-8\epsilon_i)},$
\item $\prod_{i=N_0}^\infty (1-4\epsilon_i) < \varepsilon $,
\item $\mu_{(\bC)^{\Ff}}(p)\cdot
{\prod_{i=N_0}^\infty (1-8\epsilon_i)} - \sum_{i=N_0}^\infty
\epsilon_i \geq (1-\varepsilon)\cdot \mu_{(\bC)^{\Ff}}(p),$
\item $h(N)>0$ for $N\geq N_0$.
\end{enumerate}
For $N \in \omega $  let $m^N$ be the distribution on $2^{I_0 
  \cup \dots \cup I_{N-1}}$ defined as 
$$m^N(s) = \left\{
  \begin{array}{ll}
2^{-|\bigcup_{i <
  N} I_i|}& \text{if }s \in p^N\\
0& \text{otherwise}
  \end{array} \right. .$$
Note that $\wh{m^N}$ is the counting measure of $p^N$.

Use \ref{crucialgeneralmore} to 
find $\bF^\star_N \in \FF^{N_0,N}_{\bF\rest N, h\rest N - \sr(h\rest N)}$ and
$U^\star_N \subseteq 2^{I_0 \cup \dots \cup I_{N-1}}$ such that 
\begin{multline*}\wh{\left(\MU{U^\star_N}^N\right)_{(\bC)^{\bF^\star_N}}} =
\frac{\left|p^N \cap U^\star_N \cap
    (\bC)^{\bF^\star_N}\right|}{\left| (\bC)^{\bF^\star_N}\right|}
 \geq \frac{\left|p^N \cap
    (\bC)^{\bF\rest N}\right|}{\left| (\bC)^{\bF\rest N}\right|}\cdot
{\prod_{i=N_0}^\infty (1-8\epsilon_i)} - \sum_{i=N_0}^\infty \epsilon_i,
\end{multline*}
and for  $M_0 \in
[N_0,N)$, $s \in p_s^N \rest I_0\cup \dots \cup I_{M_0-1}$ and $\bG
\in \FF^{M_0,N}_{\bF^\star_N, \sr(h\rest N)}$, 
$$ \frac{\left|p_s^N \cap U^\star_N \cap (\bC)^{\bG} \right|}
{\left|(\bC)^{\bG}\right|}
 \geq \frac{\left|p_s^N \cap U^\star_N \cap (\bC)^{\bF^\star} \right|}
{\left|(\bC)^{\bF^\star}\right|}\cdot {\prod_{i=M_0}^\infty (1-4\epsilon_i)}.$$
By compactness, there exist $\bF^\star \in \FF^\omega $ and $U^\star
\subseteq 2^{<\omega}$ such that 
$$ \forall N\  \exists M \geq N \ \lft2(\bF^\star \rest N =
\bF^\star_M \rest N \ \&\ (U^\star)^N =
(U^\star_M)^N\rgt2).$$

Put $p^\star = p \cap U^\star$ and
note that, by \ref{crucialgeneralmore}, for every $N \geq N_0$ there
exists $M \geq N$ such that  
\begin{multline*}
 \frac{\left|(p^\star)^N \cap
    (\bC)^{\bF^\star\rest N} \right|}{|(\bC)^{\bF^\star\rest N}|} = 
\frac{\left|(p^M \cap U^\star_M )^N \cap
    (\bC)^{\bF^\star\rest N} \right|}{|(\bC)^{\bF^\star\rest N}|} = 
\frac{\left|(p^M \cap U^\star_M)^N\cap
    (\bC)^{\bF_M^\star\rest N} \right|}{|(\bC)^{\bF_M^\star\rest
    N}|}\geq\\
\frac{\left|p^M \cap U^\star_M \cap
    (\bC)^{\bF^\star_M} \right|}{|(\bC)^{\bF^\star_M}|} \geq
 \frac{\left|p^M \cap
    (\bC)^{\bF\rest M}\right|}{\left| (\bC)^{\bF\rest M}\right|}\cdot
{\prod_{i=N_0}^\infty (1-8\epsilon_i)} - \sum_{i=N_0}^\infty
\epsilon_i \geq\\
 \mu_{(\bC)^{\bF}}(p)\cdot
{\prod_{i=N_0}^\infty (1-8\epsilon_i)} - \sum_{i=N_0}^\infty 
\epsilon_i \geq (1-\varepsilon)\cdot\mu_{(\bC)^{\bF}}(p).
\end{multline*}
It follows that
$$\mu_{(\bC)^{\bF^\star}}(p^\star)=\lim_{N \rightarrow \infty}
\frac{\left|(p^\star)^N \cap 
    (\bC)^{\bF^\star\rest N} \right|}{|(\bC)^{\bF^\star\rest N}|} \geq 
(1-\varepsilon)\cdot\mu_{(\bC)^{\bF}}(p).$$

Suppose that $s \in (p^\star)^{M_0}$ for some $M_0\geq N_0$. As above, for $N\geq M_0$ and
$\bG \in \FF^{M_0,N}_{\bF^\star_N, \sr(h\rest N)}$, the inequality
$$ \frac{\left|p_s^N \cap U^\star_N \cap (\bC)^{\bG} \right|}
{\left|(\bC)^{\bG}\right|}
 \geq \frac{\left|p_s^N \cap U^\star_N \cap (\bC)^{\bF^\star} \right|}
{\left|(\bC)^{\bF^\star}\right|}\cdot {\prod_{i=M_0}^\infty (1-4\epsilon_i)},$$
translates to 
\begin{multline*}\forall \bG \in \FF^{M_0, \omega}_{\bF^\star, \sr(h)} \
\mu_{(\bC)^{\bG}}(p^\star_s) \geq \mu_{(\bC)^{\Ff^\star}}(p^\star_s) 
\cdot \prod_{i=M_0}^\infty (1-4\epsilon_i)\geq \\
 \frac{1}{|(C_{M_0})^{\bF^\star(M_0)}|}
\cdot \prod_{i=M_0}^\infty (1-4\epsilon_i)>0.\end{multline*}
It follows that
if $s \in (p^\star)^{M_0}$, $M_0 \geq N_0$ then
for all $ \bG \in \FF^{M_0,
  \omega}_{\bF^\star, \sr(h)} $,
$$ \mu_{(\bC)^{\bG}}(p^\star_s) \geq
(1-\varepsilon)\cdot\mu_{(\bC)^{\Ff^\star}}(p^\star_s)>0.$$ 
Define
$$h^\star(N)=\left\{
  \begin{array}{ll}
\overline{\sr(h)}(N)& \text{if } N \geq N_0\\
0& \text{otherwise}
  \end{array}\right. \text{ for } N \in \omega .$$
Suppose that $s \in (p^\star)^N$. If $N \geq N_0$ then the above
estimates show that
$$\bv{p^\star_s,\bF^\star,h^\star}_N\geq
(1-\varepsilon)\cdot\mu_{(\bC)^{\Ff^\star}}(p^\star_s)>0.$$
If $N<N_0$ then by \ref{easy2}(3),
\begin{multline*}\bv{p^\star_s,\bF^\star,h^\star}_N \geq \sum_{s \subseteq t \in
  (p^\star)^{N_0}} \bv{p^\star_t,\bF^\star,h^\star}_N = \\
\frac{w_{N}(\bF^\star)}{w_{N_0}(\bF^\star)}\cdot\sum_{s \subseteq t \in
  (p^\star)^{N_0}} \bv{p^\star_t,\bF^\star,h^\star}_{N_0}
 >0.\end{multline*}
Finally note that for $\bG \in \FF^{\omega}_{\bF^\star, h^\star}$,
\begin{multline*} \mu_{(\bC)^{\bG}}(p^\star) = \sum_{t \in
  (p^\star)^{N_0}} \mu_{(\bC)^{\bG}}(p^\star_s) \geq \sum_{t \in
  (p^\star)^{N_0}} \mu_{(\bC)^{\bF^\star}}(p^\star_s) \cdot
(1-\varepsilon) =\\
(1-\varepsilon)\cdot \mu_{(\bC)^{\bF^\star}}(p^\star) \geq
(1-\varepsilon)^2\cdot \mu_{(\bC)^{\bF}}(p) \geq (1-2 \varepsilon)\cdot\mu_{(\bC)^{\bF}}(p).
\end{multline*}
It follows that
$$\bv{p^\star,\bF^\star,h^\star} \geq (1-2 \varepsilon)\cdot\mu_{(\bC)^{\bF}}(p).~\QED$$

\begin{theorem}\label{main1}
  Suppose that   $M_0 \in \omega $, $ \varepsilon<1$  and $\mu_{(\bC)^{\Ff}}(A)=1$.
Let $p \subseteq 2^\omega $ and $h \in {\mathcal X} $ be such that 
$$\forall N \ \forall s \in (p)^{N} \ \bv{p,\bF,h}_N>0.$$

There exist $p^\star$, $h^\star \in {\mathcal X}$ and $\Ff^\star \in
\FF^{N_0, \omega}_{\bF, h-h^\star}$ such that 
\begin{enumerate}
\item $p^\star \subseteq p \cap A$,
\item $h^\star \rest M_0=h\rest M_0$,
\item $\forall N\geq M_0 \ \log_{\sr}(h^\star)(N) = \log_{\sr}(h)(N)-1$,
\item $\forall s \in p^\star \ \bv{p^\star_s,\bF^\star,h^\star}_N>0$,
\item $\forall s \in (p)^{M_0} \ \bv{p^\star_s,\bF^\star,h^\star}_{M_0} \geq 
  (1-4\varepsilon)\cdot \bv{p_s,\bF,h}_{M_0}$.
\end{enumerate}
\end{theorem}
\Proof 
Let $\alpha = \min\left\{\bv{p_s,\bF,h}_{M_0}: s \in (p)^{M_0}\right\}.$
Fix $ \varepsilon >0$ and for every $s \in (p)^{M_0}$ find $N_0^s\geq M_0$ as in \ref{main}.
Let $N_0\geq \max\left\{N_0^s: s \in (p)^{M_0}\right\}$ be such that 
$\log_{\sr}(h)(N_0)>0$. 

Fix an enumeration $\{s_i: 0<i \leq \ell\}$ of $(p)^{M_0}$, and define
sequences
$\{\bF_i,h_i : i \leq \ell\}$ and $\{p^\star_i: 0<i \leq \ell\}$ such that 
\begin{enumerate}
\item $\bF_0=\bF$, $h_0=h$,
\item $h_i \in {\mathcal X} $ for $i \leq \ell$,
\item $p^\star_i \subseteq p_{s_i} \cap A$,
\item $\Ff_{i+1} \in \FF^{{N_0},\omega}_{\Ff_i,
    h_{i}-\sr(h_i)}$,
\item $h_{i+1}(N)=\sr(h_i)(N)$ for $N \geq N_0$, $i < \ell$,
\item $\forall i \leq \ell \ \forall N<N_0 \ h_i(N)=0$,
\item $\bv{p^\star_i,\bF_i,h_i}_{M_0} \geq
(1-4\varepsilon)\cdot\mu_{(\bC)^{\Ff_i}}(p_{s_i})$,
\item $\forall N \ \forall s \in (p^\star_i)^{N} \ \bv{p^\star_i,\bF_i,h_i}_N>0.$
\end{enumerate}

Suppose that
 $\bF_i^\star$, $ h_i^\star$ are given 
for some $i<\ell$. 
Find $q_{i+1} \subseteq p_{s_{i+1}} \cap A$ such that
$\mu_{(\bC)^{\bF_i}}(q_{i+1}) \geq (1-\varepsilon) \mu_{(\bC)^{\bF_i}}(p_{s_i}).$
Let $p_{i+1}$, $\bF_{i+1}$ and $ h_{i+1}$ be obtained by
applying \ref{main} to $q_{i+1},\ 
\bF_i$ and $h_i$. 
After $\ell$ steps we have constructed functions $\bF_\ell$, $h_\ell$
and a set $p^\star=\bigcup_{i \leq \ell} p_i$.
Functions $\bF_\ell$ and $\overline{h_\ell}=h^-$ will define walues of $\bF^\star$
and $h^\star$ for $N \geq N_0$. 

Define for $N \in \omega $,
$$h^\star(N)=\left\{
  \begin{array}{ll}
h(N)&\text{if } N<M_0\\
h^-(N)& \text{if } M_0\leq N\\
  \end{array}\right. $$
and $\bF^\star(N)=\bF_\ell(N)$ for $N\geq N_0$.
It remains to define the values of
$\bF^\star(N)$ for $N<N_0$. 

Define $\bF^\star \rest N_0$  by the following requirements:
\begin{enumerate}
\item $\bF^\star \rest M_0=\bF \rest M_0$,
\item $\bF^\star \in \FF^{M_0,\omega}_{\bF, h-h^-}$,
\item for $N<N_0$ and $s \in (p^\star)^N$,
$$p^\star_s
    \cap (\bC)^{\bF^\star\rest N_0}
    \neq \emptyset \rightarrow \left(\forall \bG \in \FF^{N,
        N_0}_{\bF^\star\rest N_0, h^\star\rest N_0} \
p^\star_s \cap (\bC)^{\bG}\neq \emptyset\right).$$
\end{enumerate}

More precisely, by induction on $N \in [M_0, N_0)$ define sequences
 $\left\{\bF^{N}_i:~i~\leq~v_{N}\right\}$ and $\left\{h^{N}_i: i
   \leq v_N\right\}$  such that  
\begin{enumerate}
\item $h_0^{M_0}=h\rest N_0$, $\bF^{M_0}_0= \bF\rest N_0$, 
$\bF^{N+1}_{0}=\bF^N_{v_N}$ and $h^{N+1}_0=h^N_{v_N}$ for $N \geq M_0$,
\item $\forall N<N_0 \ \forall i \leq v_N \ h^{N}_i \rest
  N=h^{N}_0 \rest N$,
\item $h^{N}_{i+1}=h^{N}_0 \rest N{}^\frown \sr\lft1(h_i\rest
  [N,N_0)\rgt1) $ for $ i \leq v_{N}$,
\item $\bF^{N}_{i+1} \in \FF^{N,N_0}_{\bF^{N}_{i},h_i^{N}-h_{i+1}^{N}}$,
\item if $s$ is the $i$'th element of
  $(p)^{N}$ then exactly one of the following two cases holds:
  \begin{enumerate}
  \item $ \forall \bG \in \FF^{N,N_0}_{\bF^{N}_{i},h^{N}_{i}} \ (\bC)^{\bG} \cap 
  (p_{s})^{N_0} \neq \emptyset$, 
  \item $(\bC)^{\bF^{N}_{i}} \cap 
  (p_{s})^{N_0} = \emptyset$.
  \end{enumerate}
\end{enumerate}

The construction is straightforward. If case (5a) holds, then we define
$\bF^N_{i+1}=\bF^N_{i}$, otherwise there exists $\bG \in 
\FF^{N,N_0}_{\bF^{N}_{i},h_{i}^{N}}$ such that $(\bC)^{\bG} \cap 
  (p_{s})^{N_0} = \emptyset$,
and we put 
$\bF^{N}_{i+1}=\bG$.

Observe that for $N \geq M_0$,
$h^\star(N)=h^-(N)=\sr^{(l_N)}(h)(N)=h^{N+1}_{v_N}(N)$.
Therefore
we can carry out this construction provided that
$\log_{\sr}(h)(N)>0$. However, by the choice of $ {\mathcal X} $, if
$\log_{\sr}(h)(N)=0$ then $h(N)=0$ and the required condition is
automatically met.

Finally let
$$\bF^\star(N)=\left\{
  \begin{array}{ll}
\bF(N)&\text{if } N<M_0\\
\bF^{N}(N)& \text{if } M_0\leq N<N_0\\
\bF_\ell(N)&\text{if } N\geq N_0
  \end{array}\right. .$$

We will show that $p^\star$, $\bF^\star$ and $h^\star$ have the
required properties.
Conditions (1)--(3) of \ref{main1} are obvious. 

To check (5) consider 
$s \in (p^\star)^{M_0}$.
By the choice of $N_0$, $p^\star$ and $\bF_\ell$ we have
\begin{multline*}
\bv{p^\star_s,\bF^\star,h^\star}_{M_0} \geq \min\left\{\bG\in
\FF^{M_0,N_0}_{\bF^\star\rest N_0, h^\star\rest N_0}:
\frac{\left|(p^\star_s)^{N_0} \cap
    (\bC)^{\bG}\right|}{\left(\bC)^{\bG}\right|}\right\}\cdot (1-4
\varepsilon ) = \\
\min\left\{\bG\in
\FF^{M_0,N_0}_{\bF^\star\rest N_0, h^\star\rest N_0}:
\frac{\left|(p_s)^{N_0} \cap
    (\bC)^{\bG}\right|}{\left(\bC)^{\bG}\right|}\right\}\cdot (1-4
\varepsilon ) \geq \\
\bv{p_s,\bF^\star,h^\star}_{M_0}\cdot (1-4 \varepsilon ) \geq 
\bv{p_s,\bF,h}_{M_0}\cdot (1-4 \varepsilon ).
\end{multline*}

To verify (4) we have to show that
$\bv{p^\star_s,\bF^\star,h^\star}_N>0$ for $s \in (p^\star)^N$. 
If $N \geq N_0$ it follows from the construction of $\bF_\ell$. If
$N<N_0$ then
$$\bv{p^\star_s,\bF^\star,h^\star}_N \geq (1-4 \varepsilon )\cdot
\min\left\{\frac{\left|(p^\star_s)^{N_0} \cap
      (\bC)^{\bG}\right|}{\left|(\bC)^{\bG}\right|}: \bG \in
  \FF^{M_0,N_0}_{\bF^\star\rest N_0, 
    h^\star\rest N_0}\right\}.$$
By the choice of $\bF^\star\rest N_0$, for all $\bG \in
  \FF^{M_0,N_0}_{\bF^\star\rest N_0, 
    h^\star\rest N_0}$,
$$\frac{\left|(p^\star_s)^{N_0} \cap
      (\bC)^{\bG}\right|}{\left|(\bC)^{\bG}\right|} \neq 0.$$
It follows that
$\bv{p^\star_s,\bF^\star,h^\star}_N>0$.~$\QED$

\section{Definition of $ {\mathcal P} $}\label{seven}
In this section we will define a partial order ${\mathcal P} $ having
properties A0 -- A2 from \ref{newtrick}. This will conclude the proof
of \ref{biggie}. 

We start by defining a partial ordering  $ {\mathcal Q} $
that will be used in the definition of $ {\mathcal P} $.

\begin{definition}
Let $ {\mathcal Q} $ be the following partial order:

$(p,\Ff, h) \in {\mathcal Q} $ if
\begin{enumerate}
\item $p \in \perf$, $\Ff \in \FF^\omega$, $h \in {\mathcal X}  $,
\item $\left|\dom\lft1(\Ff(k)\rgt1)\right|+h(k)\leq m_k$ for every
  $k$,
\item $p \subseteq (\bC)^{\bF}$,
\item $\forall s \in p^N \ \bv{p_s,\bF,h}_N>0$.
\end{enumerate}

For $(p^1,\Ff_1, h_1), (p^2,\Ff_2, h_2) \in {\mathcal Q}$ define
$(p^1,\Ff_1, h_1) \geq (p^2,\Ff_2, h_2)$ if
\begin{enumerate}
\item $p^1 \subseteq p^2$,
\item $\Ff_1 \in \FF^\omega_{\Ff_2,h_2-h_1}$.
\end{enumerate}
  
\end{definition}

To see that $ {\mathcal Q} $ has the fusion property we define $\geq_n$:
\begin{definition}
For $ n>0$ define $(p^1,\Ff_1, h_1) \geq_n (p^2,\Ff_2, h_2)$
if
\begin{enumerate}
\item $(p^1,\Ff_1, h_1) \geq (p^2,\Ff_2, h_2)$,
\item $\forall s \in (p^2)^{n^\star}\ 
\bv{p^1_s, \Ff_1, h_1}_{n^\star} \geq (1-2^{-n-1})\cdot \bv{p^2_s, \Ff_2,
    h_2}_{n^\star}$,
\item $h_1 \rest n^\star = h_2 \rest n^\star$,
\item $\Ff_1 \rest n^\star = \Ff_2 \rest n^\star$, 
\end{enumerate}
where
$n^\star = i_{h_1}(n).$
 \end{definition}
Note that (2) implies that $(p^1)^{n^\star}=(p^2)^{n^\star}$.

\begin{lemma}\label{qfp}
  ${\mathcal Q} $ has the fusion property.
\end{lemma}
\Proof
Suppose that $\left\{(p^k,\Ff_k, h_k): k \in \omega \right\}$ is a sequence of
conditions such that 
$(p^{k+1},\Ff_{k+1}, h_{k+1})\geq_{k+1}(p^k,\Ff_k, h_k)$ for each $k$.
Let 
$n^\star(k)=i_{h_{k+1}}(k)$.
Note that $\lim_{k \rightarrow \infty} n^\star(k)=\infty$.
Define
\begin{enumerate}
\item $h=\bigcup_{k \in \omega } h_k \rest n^\star(k)$,
\item $\Ff= \bigcup_{k \in \omega } \Ff_k \rest n^\star(k)$,
\item $p= \bigcup_{k \in \omega } (p^k)^{n^\star(k)}.$
\end{enumerate}
Observe that $h$, $\Ff$ and $p$ are well defined.

Suppose that $s \in p^{n^\star(k_0)}$, $\bG \in \FF^{N,\omega}_{\bF,h}$ and $k
\geq k_0 $, and
note that
$$\frac{\left|(p_s)^{n^\star(k)} 
   \cap 
  (\bC)^{\bG\rest n^\star(k)}\right|}{\left|(\bC)^{\bG\rest
    n^\star(k)}\right|}=
\dfrac{\left|(p_s^k)^{n^\star(k)} 
   \cap 
  (\bC)^{\bG\rest n^\star(k)}\right|}{\left|(\bC)^{\bG\rest
    n^\star(k)}\right|}
\geq \bv{p^k_s,\bF_k,h_k}.$$
Therefore
$\mu_{(\bC)^{\bG}}(p_s) \geq \inf_{k}\bv{p^k_s,\bF_k,h_k}$.
Hence, 
\begin{multline*}
\bv{p_s,\bF,h}_{n^\star(k_0)} \geq \bv{p^{k_0}_s,\bF_{k_0},h_{k_0}}_{n^\star(k_0)}\cdot
\prod_{k> k_0}\left(1-\frac{1}{2^{k+1}}\right) \geq\\
\left(1-\frac{1}{2^{k_0+1}}\right)\cdot
\bv{p^{k_0}_s,\bF_{k_0},h_{k_0}}_{n^\star(k_0)}>0.
\end{multline*}
The same computation shows that 
$(p,\bF,h) \geq_k (p^k,\bF_k,h_k)$.~$\QED$


\begin{theorem}\label{main2}
  Suppose that $(p,\Ff,h) \in {\mathcal Q}$.

If $q \subseteq p$ and
  $\mu_{(\bC)^{\Ff}}(q)>0$ then  
there exist $q^\star \subseteq q$, $\bF^\star$ and $ h^\star\in
{\mathcal X} $ such that $(q^\star, \bF^\star, h^\star) \in {\mathcal Q} $ and
 $(q^\star, \bF^\star, h^\star) \geq (p,\Ff,h)$.

If $n \in \omega $ and $A \subseteq p$ is such that
$\mu_{(\bC)^{\Ff}}(A)=1$ then there exist $q^\star \subseteq p\cap
A$, $\bF^\star$ and $ h^\star\in 
{\mathcal X} $ such that $(q^\star, \bF^\star, h^\star) \in {\mathcal Q} $ and
$(q^\star, \bF^\star, h^\star) \geq_n (p,\Ff,h)$.
\end{theorem}
\Proof
The first part follows from \ref{main} and the second from
\ref{main1}.~$\QED$

The following theorem shows that $ {\mathcal Q} $ satisfies condition
A2 defined in section \ref{out}.
\begin{theorem}\label{qa2}
  For every $(p, \bF, h) \in {\mathcal Q} $, $ n \in \omega$, $X \in
  [2^\omega]^{\leq \boldsymbol\aleph_0} $, and $\t \in \perf$ such that
  $\mu(\t)>0$,
\begin{displaymath}\mu\lft2(\lft1\{z\in 2^\omega : \exists  (q, \bG,
  f) \geq_n   (p, \bF, h)\ 
 X \cup (q + 
\rationals) 
\subseteq \t+\rationals+z\rgt1\}\rgt2)=1.\end{displaymath}
\end{theorem}
\Proof 
Suppose  that $(p,\Ff,h) \in {\mathcal Q} $ and $\t$ is a perfect set of
positive measure.

We will  need  the following observation:
\begin{lemma}
$$\mu\left(\left\{z\in 2^\omega : \mu_{(\bC)^{\Ff}}\lft1(p \cap
    (\t+z)\rgt1)>0\right\}\right)>0.$$   
\end{lemma}
\Proof
Consider the space $ p \times 2^\omega $ equipped with the
product measure $(\mu_{(\bC)^{\Ff}}\rest p)~\times~\mu$.
Let $Z = \{(x,z)\in p\times 2^\omega : z \in \t+x\}$. 
Note that $\mu\lft1((Z)_x\rgt1)=\mu(\t+x)=\mu(\t)>0$ for each $x$.
By the Fubini theorem
$$\left\{z: \mu_{(\bC)^{\Ff}}\lft1((Z)^z\rgt1)>0\right\}$$
has positive measure.  
But
$$(Z)^z = \{x \in p: z \in \t+x\}=\{x\in p: x \in \t+z\} = p \cap (\t+z).~\QED$$

Let $X  
\subseteq  2^\omega $ be a countable set. Put
$Z_X = \{z \in 2^\omega : X \subseteq \t+\rationals+z\}$. Note that 
$Z_X$ has measure one. Thus, without loss of generality, we can assume that
$X= \emptyset$.

For each $s \in p$ let 
$$Z_s = \left\{z \in 2^\omega : \mu_{(\bC)^{\Ff}}\lft1(p_s
\cap (\t+z)\rgt1)>0\right\}.$$
By the lemma, $\mu(Z_s)>0$ for each $s$.
Let $Z = \bigcap_{s \in p} (Z_s+\rationals)$. This is the measure one set
we are looking for.

Fix $z \in Z$ and $n \in \omega $. Note that 
$\mu_{(\bC)^{\bF}}(\t+\rationals+z)=1$ and apply \ref{main2}.~$\QED$

\begin{definition}
  Let $ {\mathcal P}  \subseteq {\mathcal Q} \times {\mathcal Q} $ be
the collection of elements $\lft1((p^1,\Ff_1, h), (p^2,\Ff_2, h)\rgt1)$ such
that 
\begin{enumerate}
\item $ \forall k \ \dom\lft1(\Ff_1(k)\rgt1)=\dom\lft1(\Ff_2(k)\rgt1)$,
\item $ \forall k \ \forall s \in \dom\lft1(\Ff_1(k)\rgt1) \ \lft2(\Ff_1(k)(s)=1
  \text{ or } \Ff_2(k)(s)=1\rgt2)$.

\end{enumerate}
For $\lft1((p^1,\Ff_1, h_1), (q_1,\bG_1, h_1)\rgt1),
\lft1((p^2,\Ff_2, h_2), (q_2,\bG_2, h_2)\rgt1) \in {\mathcal P} $ and $n \in \omega $
define 

$\lft1((p^1,\Ff_1, h_1), (q_1,\bG_1, h_1)\rgt1)\geq_n
\lft1((p^2,\Ff_2, h_2), (q_2,\bG_2, h_2)\rgt1)$ if

$(p^1,\Ff_1, h_1) \geq_n (p^2,\Ff_2, h_2)$ and $(q_1,\bG_1, h_1) \geq_n (q_2,\bG_2, h_2)$.
\end{definition}
Strictly speaking, the partial order used in the proof of
\ref{newtrick} was a subset of $\perf \times \perf$ while $ {\mathcal
  P} $ defined above has more complicated structure. Nevertheless it
is easy to see that it makes no difference in the proof of
\ref{newtrick} as conditions A1 and A2 refer only to the first
coordinate of ${\mathcal P} $.

 \begin{lemma}
   $ {\mathcal P} $ has the fusion property.
 \end{lemma}
\Proof Follows immediately from the definition of $ {\mathcal P} $ and 
\ref{qfp}.~$\QED$

Next we show that ${\mathcal P} $ satisfies A1.
 \begin{lemma}
   For every $ \p\in
   {\mathcal P} $, $n \in \omega $ and $z \in 2^\omega $ there exists
  $\q\geq_n
\p$ such that 
$q_1 \subseteq H+z$ or $q_2 \subseteq H+z$.
 \end{lemma}
\Proof
Suppose that $\lft1((p^1,\Ff_1, h), (p^2,\Ff_2, h)\rgt1) \in {\mathcal P} $
and $z \in 2^\omega $.

{\sc Case 1}. There exist infinitely many $k$ such that $z \rest I_k
\in \dom\lft1(\Ff_1(k)\rgt1)$. 

It follows from the definition of $ {\mathcal P} $ that in this case
there exists $i \in \{1,2\}$ and infinitely many $k$ such that 
$\Ff_i(k)(z \rest I_k)=1$. In particular, since $p^i \subseteq
(\bC)^{\Ff_i}$, for every $x \in p^i$,
$$ \exists^\infty k \ x \rest I_k \not\in C_k+z \rest I_k.$$
Thus, $p^i \subseteq H+z$.

\bigskip

{\sc Case 2}. $z \rest I_k \in \dom\lft1(\Ff_1(k)\rgt1)$ for finitely many $k$.

Let $n^\star =i_h(n)$.
Define for $ k \in \omega $, and $i=1,2$
$$\bG_i(k)=\left\{
  \begin{array}{ll}
\Ff_i(k)& \text{if } k \leq n^\star\\
\Ff_i(k) \cup (z \rest I_k, 0) & \text{if } k > n^\star
  \end{array}\right. ,$$ 
$q_i = p^i \cap (\bC)^{\bG_i}$
and
$$f(k)=\left\{
  \begin{array}{ll}
h(k)& \text{if } k \leq n^\star\\
\sr\lft1(h(k),k\rgt1) & \text{if } k > n^\star
  \end{array}\right. .$$

Clearly $\lft1((q_1,\bG_1, f), (q_2,\bG_2, f)\rgt1)\geq_n
\lft1((p^1,\Ff_1, h), (p^2,\Ff_2, h)\rgt1)$ and the same argument as in the
first case shows that it has the required properties.~$\QED$

Next we show that   $ {\mathcal P} $ satisfies A2.

 \begin{theorem} 
For every $\p\in {\mathcal P} $, $ n \in \omega
  $, $X \in
  [2^\omega]^{\leq \boldsymbol\aleph_0} $, $i=1,2$ and $\t \in \perf$ such that
  $\mu(\t)>0$,
\begin{displaymath}\mu\lft2(\lft1\{z\in 2^\omega : \exists \q \geq_n \p\  
 X \cup (q_i + 
\rationals) 
\subseteq \t+\rationals+z\rgt1\}\rgt2)=1.\end{displaymath}
 \end{theorem}
\Proof 
Suppose that $\lft1((p^1, \Ff_1, h), (p^2, \Ff_2, h)\rgt1) \in {\mathcal P}
$, $n \in \omega $, $X
\subseteq 2^\omega $ is a countable set, and $\t$ is a perfect set of
positive measure. Without loss of generality we can assume that $i=1$.
Consider the set 
\begin{displaymath}Z=\lft1\{z\in 2^\omega : \exists  (q, \bG,
  f) \geq_n   (p^1, \bF_1, h)\ 
 X \cup (q + 
\rationals) 
\subseteq \t+\rationals+z\rgt1\}.\end{displaymath}
By \ref{qa2}, $\mu(Z)=1$. Fix $z \in Z$ and let $(p',\Ff'_1,h') \geq_n (p^1, \Ff_1,h)$
be such that $p' + \rationals \subseteq \t+\rationals+z$.
Now define $\Ff_2'$ by putting
$\Ff_2'(s)=1$ for every $s \in \dom(\Ff_1') \setminus \dom(\Ff_2)$.
Clearly,
$\lft1((q, \Ff_1', h'), (p^2, \Ff_2', h')\rgt1)$ is the condition we are
looking for.~$\QED$

{\bf Acknowledgements}
We are grateful to Andrzej Ros{\l}anowski for devoting many hours to
the much needed proofreading of this paper. Thanks to his
perseverance, hopefully the process of reading this paper does not resemble the
process of writing it.


\end{document}